\documentclass[]{article}

\usepackage[english]{babel}
\usepackage{mathptmx}
\usepackage{amsmath}
\usepackage{amsthm}
\usepackage{graphicx}
\usepackage{subcaption} 
\usepackage{enumitem}

\usepackage[margin=1.5in]{geometry}
\usepackage[labelfont={bf,small}, labelsep=space]{caption}
\usepackage[colorlinks=true, linkcolor=blue]{hyperref}

\usepackage{amssymb}
\newtheorem{definition}{Definition}
\newtheorem{theorem}{Theorem}
\newcommand{\R}{\mathbb{R}}

\usepackage[maxnames=10, minnames=10, giveninits=true, style=numeric, sorting=none]{biblatex}  
\DeclareNameAlias{author}{family-given} 
\title{Accelerating Deterministic Global Optimization via GPU-parallel Interval Arithmetic}

\usepackage{authblk}

\author[1]{Hongzhen Zhang\textsuperscript{$\dagger$}}
\author[4]{Tim Kerkenhoff\textsuperscript{$\dagger$}}
\author[5]{Neil Kichler}
\author[4]{Manuel Dahmen}
\author[2,3,4]{Alexander Mitsos}
\author[5]{Uwe Naumann}
\author[1]{Dominik Bongartz\textsuperscript{*}}
\affil[1]{Department of Chemical Engineering, KU Leuven, 3001 Leuven, Belgium}
\affil[2]{JARA-CSD, 52056 Aachen, Germany}
\affil[3]{Process Systems Enginering (AVT.SVT), RWTH Aachen University, 52074 Aachen, Germany}
\affil[4]{Institute of Climate and Energy Systems - Energy Systems Engineering (ICE-1), Forschungszentrum J\"{u}lich GmbH, 52425 J\"{u}lich, Germany }
\affil[5]{Software and Tools for Computational Engineering (STCE), RWTH Aachen University, 52074 Aachen, Germany}
\affil[$\dagger$]{These authors contributed equally to this work.}
\affil[*]{Corresponding author: dominikbongartz@alum.mit.edu}
\date{} 

\begin{document}
	\maketitle
	\begin{abstract}
	Spatial Branch and Bound (B\&B) algorithms are widely used for solving nonconvex problems to global optimality, yet they remain computationally expensive. Though some works have been carried out to speed up B\&B via CPU parallelization, GPU parallelization is much less explored. In this work, we investigate the design of a spatial B\&B algorithm that involves an interval-based GPU-parallel lower bounding solver: The domain of each B\&B node is \textit{temporarily} partitioned into numerous subdomains, then massive GPU parallelism is leveraged to compute interval bounds of the objective function and constraints on each subdomain, using the Mean Value Form. The resulting bounds are tighter than those achieved via regular interval arithmetic without partitioning, but they remain fast to compute. We implement the method into our open-source solver MAiNGO via CUDA in two manners: wrapping all GPU tasks within one kernel function, or distributing the GPU tasks onto a CUDA graph. Numerical experiments show that using more subdomains leads to significantly tighter lower bounds and thus less B\&B iterations. Regarding wall clock time, the proposed spatial B\&B framework achieves a speedup of three orders of magnitude compared to applying interval arithmetic on the CPU without domain partitioning. Among the two implementations, the one developed with CUDA graph enables higher efficiency. Moreover, in some case studies, the proposed method delivers competitive or better performance compared to MAiNGO’s default solver which is based on McCormick relaxations. These results highlight the potential of GPU-accelerated bounding techniques to accelerate B\&B algorithms.
	
	
	\end{abstract}
	
	\newpage
	\section{Introduction}
	Spatial \textit{Branch and Bound} (B\&B) algorithms \cite{FalkBnB1969, LocatelliGO2013} are widely employed for global optimization of nonconvex \textit{nonlinear program}s (NLPs). They combine \textit{branching}, \textit{bounding}, and \textit{pruning} operations to explore the whole solution space, typically represented as a tree structure. The B\&B tree expands through branching, and is pruned via \textit{lower bound test}s where tighter bounds enable earlier elimination of non-promising nodes, thereby improving algorithm efficiency. For constrained problems, pruning also occurs when subproblems are proven infeasible. However, B\&B algorithms are known to be computationally intensive, particularly for problems of high dimensionality or with complex expressions in the objective function or constraints. These challenges often arise in real-world optimization problems, rendering many of them intractable for current B\&B implementations. 
	
	Over the past decades, significant effort has been made to speed up B\&B algorithms, including by improving heuristics used in searching \cite{SCHOLL199950} and branching \cite{Linderoth1999, Morrison2014}, constructing tighter relaxations for bounding \cite{He2020, Ye2025}, reducing the node domain via bound tightening techniques \cite{Puranik2017, schichl2005interval} or reformulating to a reduced space by eliminating variables \cite{Bongartz2017, mitsos2009mccormick}. With the emergence of computing architectures such as multi-core \textit{Central Processing Units} (CPU) and \textit{Graphics Processing Units} (GPU), one promising alternative to speed up B\&B algorithms is parallelization.
	
	CPU parallelization of B\&B algorithms has been attracting interest for several decades already \cite{eckstein1994parallel, bader2005parallel, talbi2006parallel, langer2013parallel}: One of the first papers on parallel B\&B algorithms was published by Mohan in 1982 \cite{mohan1982study}, where the traveling salesman problem, an \textit{mixed-integer programming} (MIP) problem, was discussed. Later, these efforts were extended to nonlinear and mixed-integer nonlinear programming problems \cite{Pardalos1989, Smith1999}. On modern computing architectures, significant speedups can be achieved via CPU parallelization. For example, Archibald et al. \cite{archibald2018replicable} achieved more than 100$\times$ speedup for the maximum clique problem using a Beowulf cluster with up to 200 CPU cores. Similarly, in our previous work \cite{dominik2021}, our open-source solver MAiNGO \cite{S1_bongartz_2018_maingo} achieved approximately 300$\times$ speedup on 1500 CPU cores when solving nonconvex NLPs. However, while many modern MILP solvers now support CPU parallelization, many global NLP solvers still lack this feature. Furthermore, parallel performance is often highly problem-dependent, and access to high-performance computing (HPC) clusters may be limited or unavailable to general users. 
	
	\label{para:parallelization_types}
	Most CPU parallelization schemes for B\&B algorithms can be classified into one or more of three main categories \cite{gendron1994parallel}: \textbf{1.} introduce parallelism at the "operation level", e.g., bounding multiple existing B\&B nodes in parallel to accelerate the solution of the subproblems \cite{Pekny1992}. \textbf{2.} Build the B\&B tree in parallel (node-level parallelism), i.e., each processing unit will take one B\&B node and perform all necessary operations \cite{martinez2004interval, casado2008branch}. Since both branching and pruning operations will be performed, the nodes being processed simultaneously can either further grow to subtrees (with two child nodes) or be pruned. \textbf{3.} Run multiple B\&B trees in parallel, where one tree can vary from the other by modifying the branching and bounding operations, or by applying different searching strategies \cite{janakiram1988randomized}, etc.. 
	
	Compared to CPU parallelization, GPU parallelization of B\&B is even less explored. In contrast, in other areas of scientific computing, it has already shown great promise, including recent successful efforts to accelerate local optimization methods \cite{pacaud2024gpuLocal, cole2023exploiting}. Although in principle, the aforementioned CPU parallelization strategies for B\&B can also be adapted for GPUs, the inherent irregularity of the B\&B search tree poses significant challenges for efficient implementation on the \textit{single instruction multiple threads} (SIMT) architecture of modern GPUs. 
	
	The first B\&B algorithms leveraging GPU parallelization were proposed in 2013 by Kiel et al. \cite{kiel2013uses} and Eriksen et al. \cite{eriksen2013gpu}, and were implemented with the GPU programming framework CUDA \cite{CUDA}. Both of them were designed for parameter estimation, an unconstrained NLP. In 2016, Vu and Derbel \cite{vu2016parallel} explained the hurdles of designing a parallel B\&B algorithm on heterogeneous HPC platforms. Although they designed the parallel B\&B algorithm for MIPs, the general framework proposed could be adapted to solving nonconvex NLPs. 
	
	\label{para:GPU_parallelization_limitation}
	A recent approach for more general NLPs by Gottlieb et al. \cite{Gottlieb2024} proposed to calculate the lower bounds of numerous B\&B nodes in parallel on the GPU, using McCormick relaxations \cite{mccormick1976computability, tsoukalas2014multivariate} via automatic source code generation. Very recently, Zhang et al. \cite{zhang2025gpu} proposed a similar method for unconstrained problems. Like Gottlieb et al. \cite{Gottlieb2024}, they also bound multiple B\&B nodes in parallel on the GPU, but with Natural Interval Extension \cite{moore1966interval} instead of McCormick relaxations. While these two approaches enabled noticeable speedups, there are some aspects that can possibly limit their performance: The McCormick relaxations evaluated on one B\&B node can differ from those evaluated on others, which implies potential imbalance and thus \textit{warp divergence} between GPU threads. Though interval computations generally involve less divergence, the bounds derived via \textit{Natural Interval Extension} can be quite loose, thus hindering the effective pruning of non-promising B\&B nodes. In addition, at the start-up phase of B\&B methods, only few B\&B nodes are generated and therefore the massive GPU parallelism cannot be fully utilized. Even though one can partition the original domain to produce a large amount of nodes at the very beginning \cite{Gottlieb2024, zhang2025gpu}, such an extensive a priori partitioning may increase the size of the B\&B tree disproportionally and thus slow down the solution. 
	
	In this work, we thus propose an alternative approach in which we do not aim to parallelize the B\&B algorithm itself, but rather to view the GPU as a tool to quickly compute tight bounds, such that less B\&B iterations are required to converge. In particular, we develop a lower bounding solver (assuming a minimization problem) that takes inspiration from the \emph{splitting technique} known from interval methods: it only \emph{temporarily} partitions the B\&B node into a large number of subdomains and evaluates interval bounds on each subdomain, to obtain an overall bound on the  B\&B node that is tighter than without partitioning. By conducting the splitting technique in parallel on the GPU, this tighter bound remains fast to compute. Instead of Natrual Interval Extension, we use the \textit{Mean Value Form} \cite{moore2009introduction} to derive the interval bounds on each subdomain. 
	
	This alternative approach offers several advantages: First, with interval methods, warp divergence is likely reduced compared to McCormick relaxations, as fewer execution branches are involved. Second, the Mean Value Form can attain higher convergence order than the Natural Interval Extension \cite{hansen1969topics} and thus is likely to benefit even more from the splitting technique in terms of bound tightness. Finally, the GPU parallelization is completely decoupled from the state of the B\&B tree, and any desired number of subdomains can be generated to saturate the available GPU cores as much as possible without partitioning the B\&B tree massively at the beginning of the algorithm. As additional improvements, we explore the use of CUDA Graphs for more efficient evaluation on the GPU. 
	
	The reminder of this paper is organized as follows: In Section 2, the notation and necessary background on interval arithmetic and GPU parallelization are introduced. In Section 3, the workflow of our method and the implementation details are discussed. In Section 4, the performance of our method is demonstrated on a number of numerical experiments. Finally, the conclusions are summarized and future research directions are discussed in section 5.
	
	\section{Background}
	In the proposed method, the bounds are obtained via interval arithmetic enhanced with GPU parallelization. In this section, after some preliminaries, we thus introduce background on interval arithmetic  including the Natural Interval Extension, the Mean Value Form, and the splitting technique underlying our lower bounding method, as well as relevant aspects of GPU parallelization and the CUDA toolkit.
	
	\subsection{Preliminaries} 
	We adopt the convention of denoting intervals by capital letters, and we only consider compact intervals. The endpoints of an interval are wrapped within square brackets, and the elements of a multi-dimensional vector are encapsulated within parenthesis. For example, a non-empty one-dimensional interval $X$ is denoted by $X = \left[ \underline{X}, \overline{X} \right]$ with $\underline{X} \leq \overline{X}$, where the underline and overline denote the left and right endpoints of an interval, respectively. An n-dimensional interval box $\textit{\textbf{X}}$, is defined by $\textit{\textbf{X}} = (X_{1}, X_{2}, ..., X_{n})$ and each $X_{i}$ is an one-dimensional interval. We denote by $\underline{\textit{\textbf{X}}} = (\underline{X_{1}}, \underline{X_{2}}, ..., \underline{X_{n}})$ the vector containing the left endpoints of $X_{i}$, similarly, $\overline{\textit{\textbf{X}}} = (\overline{X_{1}}, \overline{X_{2}}, ..., \overline{X_{n}})$ the vector containing the right endpoints of $X_{i}$. Let $\textit{\textbf{D}}\subseteq\R^n$. Then, $I(\textit{\textbf{D}})$ denotes the set of all compact interval subsets of $\textit{\textbf{D}}$. 
	
	For a one-dimensional interval $X$, the \textit{width} of $X$ is defined and denoted by $w(X) := \overline{X} - \underline{X}$ and the \textit{midpoint} of $X$ is given by $m(X) := \frac{1}{2}(\underline{X} + \overline{X})$. These notions can be extended to multi-dimensional interval boxes: The width of an interval box $\textit{\textbf{X}}$ is $w(\textit{\textbf{X}}) := \max\limits_{i} \, w(X_{i})$ and the midpoint of an interval box $\textit{\textbf{X}}$ is given by $m(\textit{\textbf{X}}) := (m(X_{1}), m(X_{2}), ..., m(X_{n})).$
	
	The \textit{intersection} and \textit{union} of two intervals are straightforward. The union of two intervals is generally a disconnected set and thus not an interval. However, the \textit{interval hull} of two intervals $$X \, \underline{\cup} \, Y := [\min\left\{\underline{X}, \, \underline{Y}\right\}, \max\left\{\overline{X}, \, \overline{Y}\right\}],$$ where the min and max operations are to be applied component-wise in the case of interval boxes, is always an interval. Clearly, $X \cup Y \subseteq X \, \underline{\cup} \, Y.$ 
	
	Considering a real-valued function $f$ defined on an n-dimensional domain $\mathbf{D}$, its image set over a subset $\textit{\textbf{X}}\subseteq\mathbf{D}$ is denoted as $f(\textit{\textbf{X}}) := \left\{ f(\textit{\textbf{x}}): \textit{\textbf{x}} \in \textit{\textbf{X}}\right\}.$
	
	\subsection{Interval Arithmetic}
	\label{subsec:intervals}
	The concept of interval arithmetic, which replaces real-valued operations with suitable interval-valued operations, was originally proposed to bound rounding errors in floating-point computations \cite{moore1966interval}. Hansen \cite{hansen1979global, hansen1980global} extended it to a computationally cheap lower bounding method for global optimization.
	
	In this context, bounds on a function $f$ are obtained through an inclusion function:
	
	\begin{definition}[Inclusion function \cite{LocatelliGO2013}]
		Given $f: \textit{\textbf{D}}\subseteq \R^n\rightarrow \R$, a function  $F:I(\textit{\textbf{D}})\rightarrow I(\R)$ if called an 
		\emph{inclusion function of $f$} if for any $\textit{\textbf{X}} \in I(\textit{\textbf{D}})$, $ 
		f(\textit{\textbf{X}})  \subseteq F(\textit{\textbf{X}})$.
	\end{definition}
	However, for an inclusion function $F$ to be useful in B\&B, it needs to converge to $f$ as $\textit{\textbf{X}}$ shrinks, which is the case for inclusion isotonic interval extensions:

	\begin{definition}[Inclusion isotonicity \cite{moore2009introduction}]
		A function $F:\textbf{D}\in I(\R^n)\rightarrow I(\R)$ is called \emph{inclusion isotonic} if for any $\textit{\textbf{X}}, \textit{\textbf{Y}} \in I(\textit{\textbf{D}})$, we have
		$\textit{\textbf{Y}} \subseteq \textit{\textbf{X}} \implies F(\textit{\textbf{Y}}) \subseteq F(\textit{\textbf{X}})$.
	\end{definition}
	
	\begin{definition}[Interval extension \cite{moore2009introduction}]
		Given $f: \textit{\textbf{D}}\subseteq \R^n\rightarrow \R$, a function $F:I(\textit{\textbf{D}})\rightarrow I(\R)$ is called an \emph{interval extension of $f$} if for any degenerate interval $[\textit{\textbf{x}}, \textit{\textbf{x}}] \in I(\textit{\textbf{D}})$,
		$F([\textit{\textbf{x}}, \textit{\textbf{x}}]) = f(\textit{\textbf{x}}).$
	\end{definition}

	Interval arithmetic allows computing such inclusion functions of factorable functions:
	\begin{definition}[Factorable function]
		A function is \emph{factorable} if it can be expressed as a finite composition of functions from a given library of so-called \emph{intrinsic functions}, for which the image set over any interval subset of their domain is known.
	\end{definition}
	Typical examples for intrinsic functions are binary operations such as addition and multiplication, or simple functions such as $\log(\cdot)$, $\sin(\cdot)$, or $\tanh(\cdot)$. For the former, the image sets are given by the corresponding \emph{interval binary operations} \cite{scholz2012theoretical}, while for the latter they are sometimes called \emph{united extension}. Functions for computing these image sets have been implemented in multiple software libraries, see \cite{althoff2016implementation} for an example.

	Given a factorable function $f$, the most straighforward way to obtain an inclusion function $F$ via interval arithmetic is the \emph{Natural Interval Extension}, where each occurrence of the independent variables is replaced with the corresponding intervals, and for each occurrence of an intrinsic function, its image set is computed over the interval of its argument. While cheap to compute, the bounds from the Natural Interval Extension are often rather weak.
	
	One alternative is the so-called \emph{Mean Value Form}, which is often regarded as a special case of a \textit{Centered Form} \cite{hansen1969topics}:
	\begin{definition}[Mean Value Form]
		Let $f:D\subseteq\R^n\rightarrow\R$ be a factorable function and $\nabla F(\textit{\textbf{X}})$ the Natural Interval Extension of $\nabla f(\textit{\textbf{x}})$ over $\textit{\textbf{X}} \in I(D)$, then the Mean Value Form of $f$ over $\textit{\textbf{X}}$ is $$F(\textit{\textbf{X}}):=f(\mathbf{m}(\textit{\textbf{X}})) + \langle\nabla F(\textit{\textbf{X}}), \textit{\textbf{X}} - \mathbf{m}(\textit{\textbf{X}})\rangle.$$
		\label{def:MVF}	
	\end{definition}
	Since $\nabla F(\textit{\textbf{X}})$ is derived via Natural Interval Extension and thus an inclusion function of $\nabla f(\mathbf{x})$, by the mean value theorem, $f(\textit{\textbf{X}}) \subseteq  F(\textit{\textbf{X}})$,
	i.e., the Mean Value Form also provides an inclusion function of $f$.

	Compared to Natural Interval Extension, the Mean Value Form tends to produce tighter interval bounds as the interval width decreases, leading to faster convergence in optimization. Hansen \cite{hansen1969topics} proved that the Mean Value Form can attain a convergence order of two, whereas Natural Interval Extension can only attain a convergence order of one. 
	
	In B\&B algorithms, the above properties that ensure the computed bounds converge to the function $f$ when constructed over smaller intervals are key for ensuring convergence. Tigther bounds usually result in less B\&B iterations until convergence. However, whether or not the wall-clock time also decreases depends on the time required to compute the bounds. In this work, we address this challenge by developing a method to compute tighter bounds through a temporary partitioning of the domain of a B\&B node, while keeping the computational time for bounding low through parallelization, as will be elaborated on below.
	
	\label{para:splitting}
	In interval methods, obtaining tighter bounds via (uniform) partitioning is known as the \emph{splitting technique} \cite{moore2009introduction}:
	\begin{enumerate}[label=\roman*:]
		\item For a function $f$ defined on an $\textit{\textbf{X}} = [X_{1}, X_{2}, ..., X_{n}]$, we can obtain $N^n$ subdomains by partitioning each dimension uniformly into $N$ subintervals. Specifically, for each dimension $X_{i}$, we get the subintervals $$X_{i,j} = [\underline{X}_{i} + (j-1) \, w(X_{i})/N, \underline{X}_{i} + j \, w(X_{i})/N],$$ where $j = 1, 2, ..., N$. Thus one subdomain $k$ can be denoted by $$\textit{\textbf{X}}_{k} = (X_{1, j_{1}}, X_{2, j_{2}}, ..., X_{n, j_{n}}).$$ Note that $j_{i} \in \left\{ 1, 2, ..., N \right\}$, i.e., each subdomain $\textit{\textbf{X}}_{k}$ is constructed by selecting one subinterval $X_{i,j}$ on each dimension $i$. Therefore we have $$\textit{\textbf{X}} = \cup_{k=1}^{N^n} \textit{\textbf{X}}_{k},$$ with $w(\textit{\textbf{X}}_{k}) = w(\textit{\textbf{X}})/N$. 
		\item The inclusion function $F$ derived via either Natural Interval Extension or Mean Value Form can be evaluated on each subdomain $\textit{\textbf{X}}_{k}$.
		\item By taking the interval hull of all evaluated $F(\textit{\textbf{X}}_{k})$, we obtain the interval bounds of $f$ over the whole domain $\textit{\textbf{X}}$, $$F_{(N)}(\textit{\textbf{X}}) = \underline{\cup}_{k=1}^{N^n} F(\textit{\textbf{X}}_{k}).$$ $F_{(N)}(\textit{\textbf{X}})$ is also called a \textit{refinement} of $F(X)$. 
	\end{enumerate}
	
	\label{para:excesswidth}
	Given an inclusion function derived via Natural Interval Extension or Mean Value Form, for which over its interval domain $\textit{\textbf{X}}$, no domain violations occur during the evaluation, and for which all intrinsic functions are Lipschitz,  the improvement in bound tightness through the splitting technique can be quantified as follows:
	\setcounter{theorem}{0}
	\begin{theorem}[Refinement excess width \cite{moore2009introduction}]
		If $F$ is an inclusion function derived from either Natural Interval Extension or Mean Value Form, then $$F_{(N)}(\textit{\textbf{X}}) = f(\textit{\textbf{X}}) + E_{N},$$ and $w(E_{N})$ is called the excess width of $F_{(N)}(\textit{\textbf{X}})$. For Natural Interval Extension, some $K \in \mathbb{R}$ exists such that $w(E_{N}) \leq K \, w(\textit{\textbf{X}})/N$. For Mean Value Form, some $\tilde{K} \in \mathbb{R}$ exists such that $w(E_{N}) \leq \tilde{K} \, (w(\textit{\textbf{X}})/N)^2$.
		\label{theorem:splitting}
	\end{theorem}
	This shows that the overall interval bounds obtained via splitting can be tighter compared to those calculated by applying interval arithmetic directly over the entire $\textit{\textbf{X}}$, which is the cornerstone of the new GPU-parallel lower bounding method to be presented in this work. The Mean Value Form is further expected to benefit more from the splitting technique because of the quadratic improvement in excess width. Note however that because in general $K\neq\tilde{K}$, it may not always produce tighter bounds compared to Natural Interval Extension when the interval width is still large.
	
	Temporary partitioning to achieve tighter lower bounds has also been used in the form of piece-wise linear relaxations, e.g., for quadratically constrained problems \cite{hasan2010piecewise,misener2011apogee}. However, these piece-wise linear relaxations result in a mixed-integer problem to compute the overall bound. In contrast, the splitting technique for interval methods merely requires taking the interval hull, which is much cheaper, thus leaving the bulk of the computation in the evaluation of the bounds over the subdomains. As such, it lends itself more naturally to the GPU parallelization scheme developed herein.

	\subsection{GPU Parallelization and CUDA Toolkit}
	In contrast to CPUs, GPUs are designed with a large number of lightweight cores which are organized into streaming multiprocessors, making them well-suited for data-parallel computations. NVIDIA GPUs, for example, execute parallel tasks on threads grouped into warps, which are scheduled across streaming multiprocessors under the SIMT architecture. SIMT helps to improve throughput and hide latency, but it can become a bottleneck when threads within a warp follow divergent execution paths, a phenomenon known as warp divergence. Consequently, algorithms originally designed for CPUs usually require adaptation to align with the SIMT architecture such that warp divergence can be reduced.
	
\begin{figure}[htbp]
		\centering
		\includegraphics[width=\textwidth, trim=0.5cm 1cm 1cm 1cm, clip ]{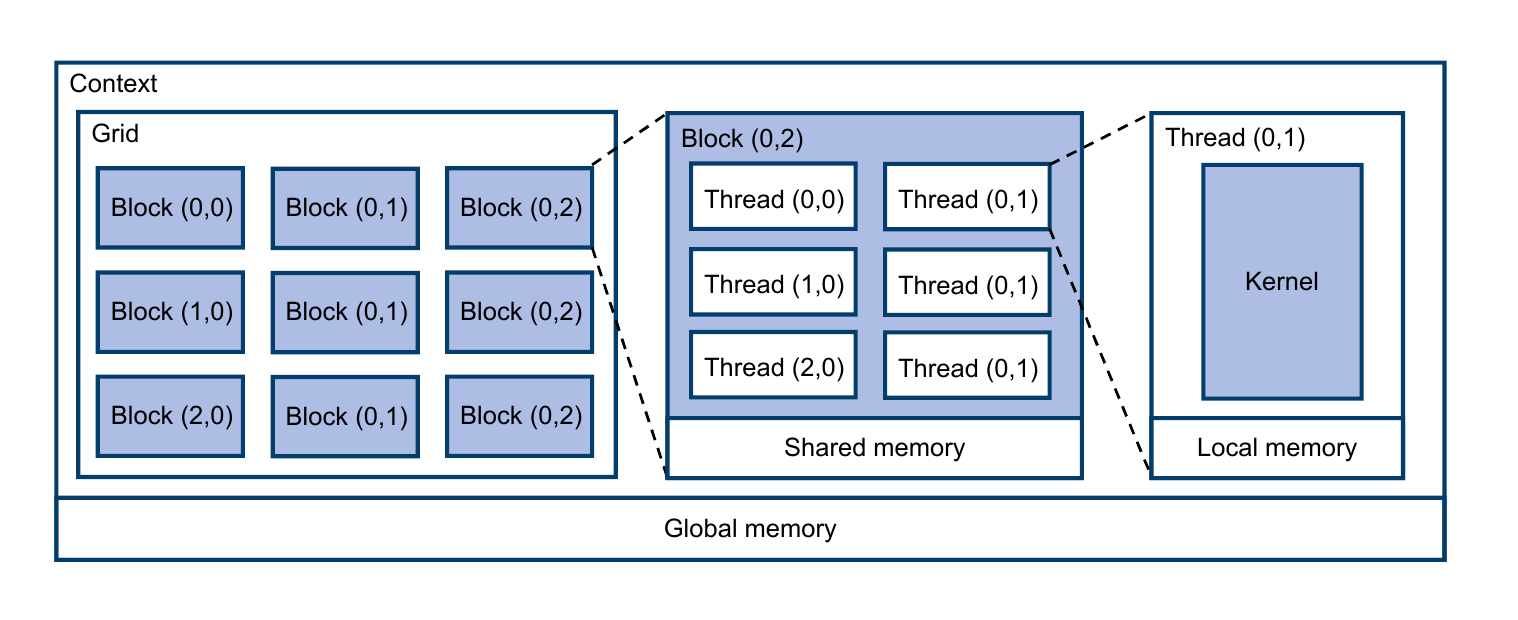}
		\caption{\footnotesize CUDA provides a programming model to allocate tasks onto GPU resources. With CUDA, one can request a certain number of CUDA cores on the GPU by specifying the number of blocks $n_{b}$ and the block size $n_{t,b}$ (number of threads per block) that are required by a given computation task. A reasonable number of threads should be requested such that they can be executed efficiently by the CUDA cores available on the GPU. Otherwise, an excessive number of threads may lead to serialized execution due to resource limitation, whereas requesting too few of them may result in a waste of GPU resources. 
		}
		\label{fig:CUDA_programming}
	\end{figure}

	For heterogeneous environments containing both CPUs and GPUs, several platforms \cite{Karimi2010} have been developed or extended to enable efficient programming. The most widely used is CUDA, a general purpose parallel computing platform and programming model. It comes with a software environment that allows developers to use C++/C as a high-level programming language. In CUDA, the function to be executed in parallel is called a \emph{kernel}, and it is mapped onto a grid consisting of multiple blocks. One block contains multiple threads, and each thread will execute a single instance of the kernel on the data it possesses, see Figure~\ref{fig:CUDA_programming}. In a heterogeneous environment, data transfers are indispensable around the execution of kernels, since CPUs and GPUs possess different memory spaces. CUDA offers a suite of APIs to manage memory allocation and data movement between different devices.
	
	Besides submitting each GPU task separately via kernels, CUDA offers an alternative model for workload submission called CUDA Graphs. They are constructed by connecting a series of operations, e.g., kernels, CPU functions, and memory operations, via dependencies. Workload submission using CUDA Graphs is divided into three distinct stages: definition, instantiation, and execution, therefore a CUDA Graph can be defined once and then launched repeatedly. 
	
\begin{figure}[htbp]
		\centering
		\includegraphics[width=0.45\linewidth, trim=6cm 0cm 6cm 0cm, clip]{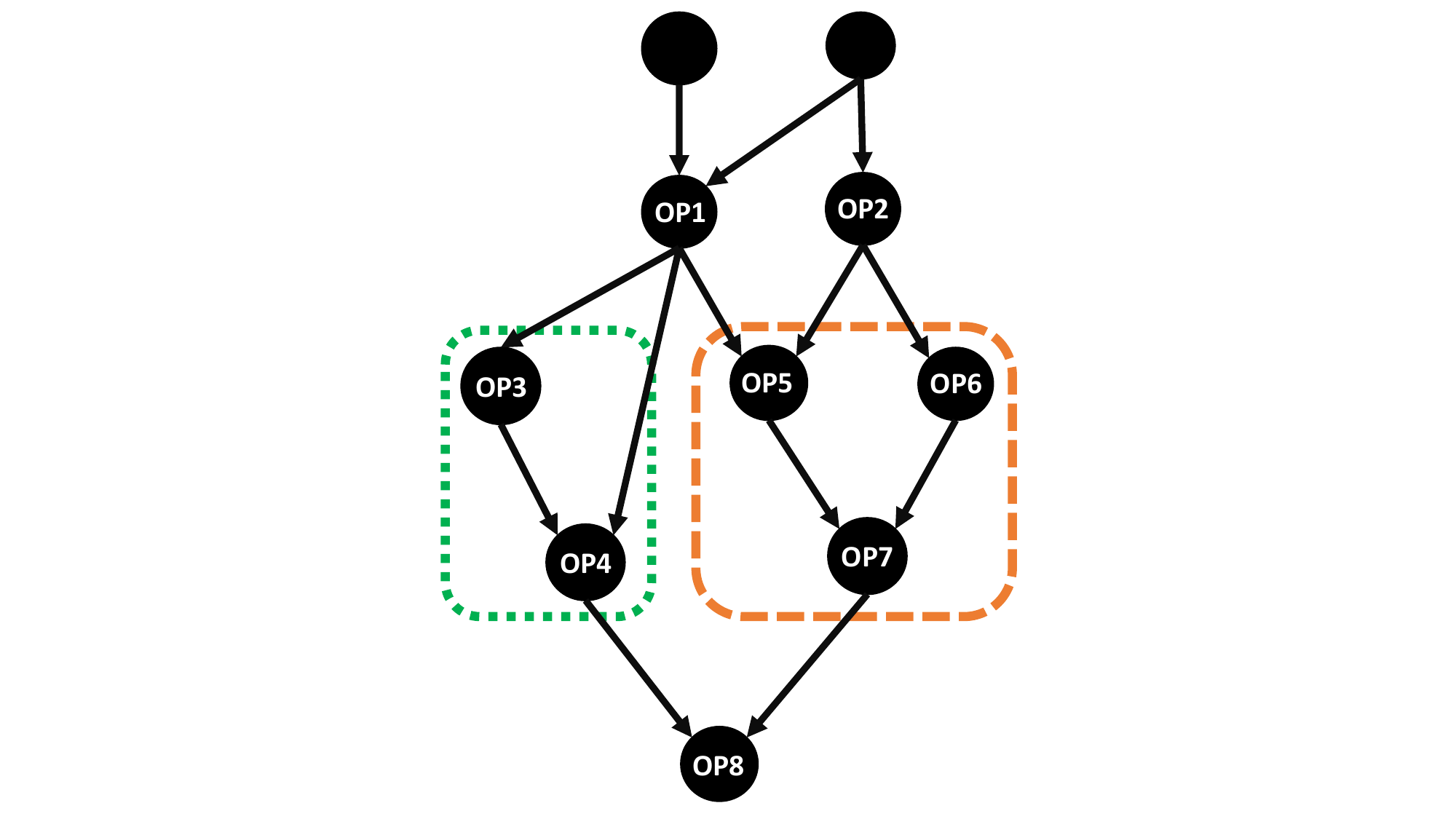}
		\caption{\footnotesize CUDA Graphs potentially offer more parallelism: In this example consisting of a CUDA Graph with 8 operations, operations OP3-OP4 in the green box and operations OP5-OP7 in the orange box are independent from each other. These two sets of operations can possibly run in parallel if the resources allow.}
		\label{fig:CUDA_graph}
	\end{figure}

	Given a computational task consisting of multiple kernel calls, CUDA Graphs offer at least two ways to improve the efficiency. Considering the execution of multiple kernels, a natural way is submitting them to the GPU in a piece-wise manner. In this case, before each submission, the host device, i.e., the CPU, has to perform a sequence of operations in preparation for setting up and launching the kernel on the GPU side, leading to significant overhead on the host (CPU) side especially when the actual execution of kernels takes only a short period of time. Contrarily, after defining a CUDA graph containing these kernels, one can submit all of them to the GPU through only one launch of this graph, thus the preparation overhead between kernel submissions can be saved and the CPU will be freed when the GPU is processing the submitted tasks. On the other hand, independent branches within the same CUDA graph can run in parallel if the resources allow, which (potentially) results in more parallelism (see Figure~\ref{fig:CUDA_graph}).
	
	GPUs can have computationl advantages over CPUs specifically in the context of interval arithmetic, where directed rounding is usually desired to guarantee inclusion in interval operations: On x86-64 CPUs, the global rounding mode has to be modified to perform directed outward rounding. This is expensive because it requires to flush the CPU instruction pipeline. In $2012$, Collange et al. \cite{collange2012interval} introduced the first interval arithmetic library for CUDA and demonstrated that this GPU interval library can be more efficient than the Boost Interval library on the CPU \cite{bronnimann2006design}. However, since the library developed by Collange et al. covers only a very limited number of intrinsic functions, an extended CUDA interval arithmetic library \cite{Kichler2025cuinterval} is implemented for this work as will be detailed later.
	
	\section{Methodology and Implementation}
	Due to the irregularity of B\&B search trees as well as the differences between CPU and GPU architectures, established CPU parallelization schemes for B\&B algorithms are generally not well suited for the SIMT architecture of GPUs. Thus, and to address the challenges of existing GPU parallelization schemes outlined above (cf. Section~\ref{para:parallelization_types}), we propose to parallelize only the lower bounding part of the B\&B algorithms on the GPU, as a way to derive tighter bounds and thus reduce the number of B\&B iterations, while keeping the time for bounding short. This also ensures that the available GPU cores can be close to fully utilized irrespective of the size of the B\&B tree. 
	
\begin{figure}[htbp]
		\centering
		\includegraphics[width=\textwidth, trim=0cm 3cm 0cm 3cm, clip]{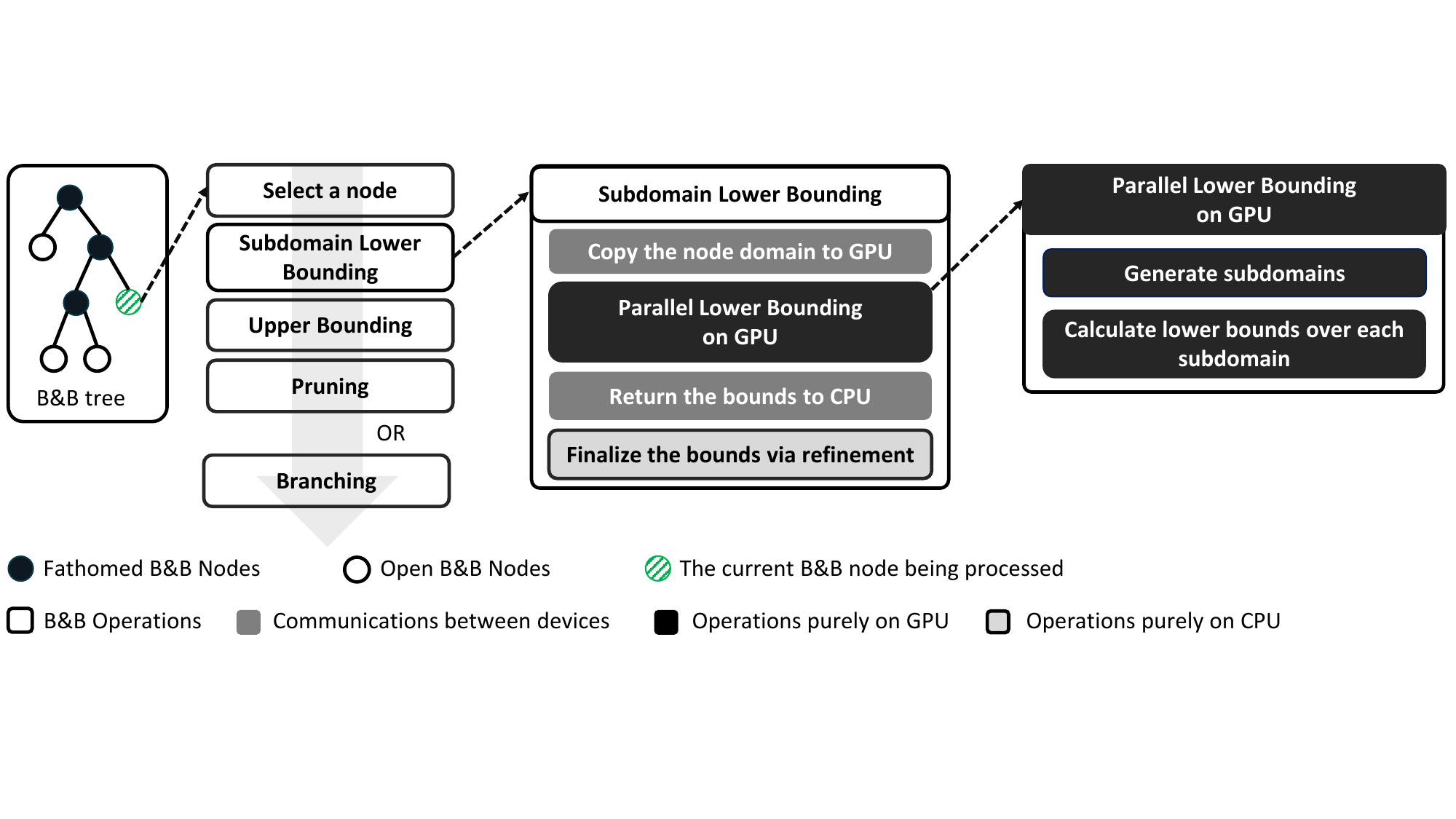}
		\caption{\footnotesize The proposed heterogeneous spatial B\&B algorithm. Rather than processing multiple B\&B nodes in parallel, we parallelize only the lower bounding solver on the GPU to derive tighter lower bounds and reduce the number of B\&B iterations, such that the parallelization does not depend on the state of the B\&B tree, and in particular does not require many open B\&B nodes being available. In the proposed Subdomain Lower Bounding method, a B\&B node is temporarily split into many subdomains, and an inclusion function is evaluated on all subdomains, then the corresponding refinement is derived. This results in much tighter bounds, but the evaluation remains fast thanks to the massive GPU parallelism.}
		\label{fig:BnB_decompose}
	\end{figure}

	In our proposal, the splitting technique (cf. Section~\ref{para:splitting}) originating from interval arithmetic is extended into an efficient GPU-parallel lower bounding method, referred to as \textit{Subdomain Lower Bounding}. At each B\&B iteration, GPU-related tasks only occur during the Subdomain Lower Bounding phase. Other operations, including node selection, upper bounding with local solver, and pruning (via either lower bound test or feasibility check) or branching, are managed by the CPU. The entire heterogeneous parallelization scheme for spatial B\&B is shown in Figure~\ref{fig:BnB_decompose}. 
	
	\subsection{The Subdomain Lower Bounding Method}
\begin{figure}[htbp]
		\centering
		\includegraphics[width=\textwidth, trim=0cm 3.5cm 0cm 4cm, clip]{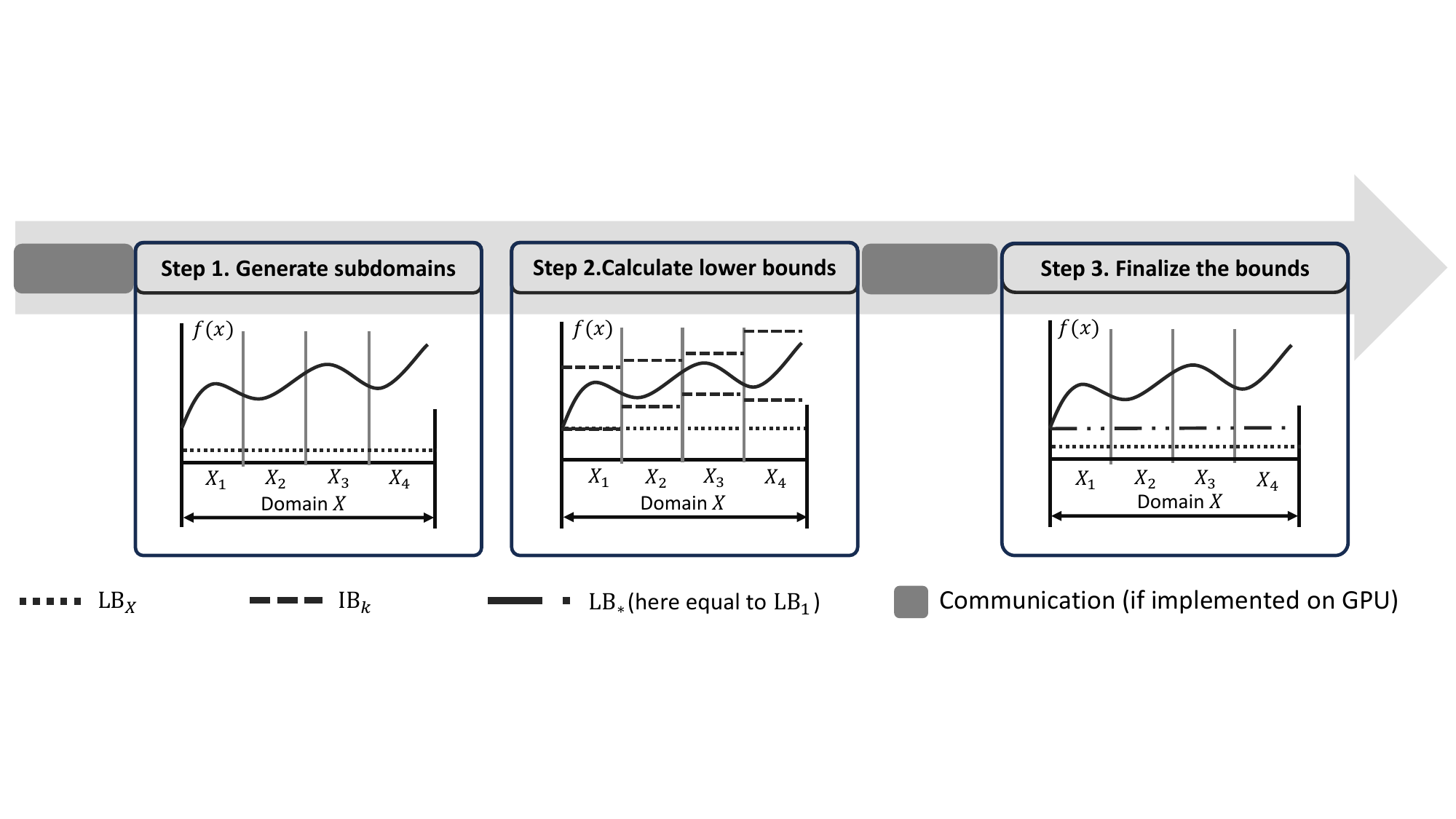}
		\caption{\footnotesize Our Subdomain Lower Bounding method consists of three main steps, illustrated here for a one-dimensional problem with only box constraints: \textbf{Step 1.} Generate subdomains by partitioning the original domain $X$ of a B\&B node into smaller subdomains $X_k$ ($k=\textup{1, 2, 3 ...}$). \textbf{Step 2.} Perform Natural Interval Extension or Mean Value Form to calculate the interval bounds $\textup{IB}_k$ of function $f(x)$, over each subdomain $X_k$. \textbf{Step 3.} Identify the refinement (the interval hull of all $\textup{IB}_k$). The left endpoint of the refinement is used as the overall lower bound $\textup{LB}_*$ of function $f$ on domain $X$. Compared to the lower bound obtained by applying interval arithmetic directly on domain $X$ ($\textup{LB}_X$; shown only for comparison), $\textup{LB}_*$ can be much tighter. For problems with general inequality constraints of the form $g(x)\leq0$ and $h(x)=0$, the bounds on $g$ and $h$ over $X$ are obtained analogously.}
		\label{fig:Subdomain_workflow}
	\end{figure}

	Consider an NLP of the form
	\begin{equation*}
		\label{eq:NLP}
		\begin{aligned}
			& \min_{\textit{\textbf{x}}\in \textit{\textbf{X}}} & & f(\textit{\textbf{x}}) \\
			& \text{s.t.} & & g_i(\textit{\textbf{x}}) \leq 0 \hspace{0.3cm} \text{ for } i=1, ..., n_g\\
			& & & h_j(\textit{\textbf{x}}) = 0 \hspace{0.3cm} \text{ for } j=1, ..., n_h,
		\end{aligned}
	\end{equation*}
	where \textit{\textbf{X}} can denote either the original domain of the problem, or the domain of a B\&B node. To compute the lower bound of such an NLP over $\textit{\textbf{X}}$, we apply our Subdomain Lower Bounding method, of which the three main steps are illustrated in Figure~\ref{fig:Subdomain_workflow} (for the case of $n_g=n_h=0$). 
	
	\label{para:steps_on_GPU}
	Steps 1 and 2 are designed to be executed in parallel on the GPU, thus prior to initiating the lower bounding process for each B\&B node, the domain associated with the selected B\&B node is transferred from CPU memory to GPU memory, ensuring accessibility for GPU threads. In Step 1, the domain of the selected B\&B node is partitioned into multiple subdomains according to a predefined \textit{partitioning strategy}, which will be discussed later. Subsequently, in Step 2, the GPU threads pick up subdomains (ideally exactly one subdomain per thread), and compute interval bounds of the objective function $f$ and constraints $g_i$ and $h_j$ over that subdomain. After completion of Step 2, the interval bounds computed over all subdomains are transferred back from GPU memory to the CPU, to facilitate the following calculation of refinements. In Step 3, which is executed on the CPU, the refinement is determined for each function via a sorting procedure, yielding the overall bounds of $f$, $g_i$ and $h_j$ on the entire node domain as left and right endpoints of the correspond refinements.
	
	\label{para:constraint_bound_pruning}
	The lower bound of the objective function derived from Step 3 can then be used on the CPU as in any B\&B algorithm to prune nodes via the lower bound check. Similarly, the bounds of the constraints can be used to prune nodes via the feasibility check: Given a B\&B node and the bounds of $g_i$ and $h_j$ on this node, an inequality constraint $g_i$ is violated if the lower interval bound of this constraint is larger than $0$. An equality constraint $h_j$ is violated if $0$ is not contained within its bounds.
	
	\subsubsection{Basic partitioning strategies}
	Assuming a hardware system with only one NVIDIA GPU, when using our heterogeneous spatial B\&B algorithm, it is recommended to partition the domain of a B\&B node into as many subdomains as possible, until all the CUDA cores on this GPU are saturated, such that the GPU parallelism can be fully leveraged. The two most basic partitioning strategies are the following:
	\begin{itemize}[label=\textbullet]
		\item \textit{Uniform partitioning}: each dimension of the domain of the selected B\&B node is partitioned into the same number of one-dimensional subintervals. Then all possible combinations of the obtained subintervals will be used as subdomains. E.g., given the node domain $X=[-3,1]\times[0,2]$ and assuming the desired number of subdomains is four, the uniform partitioning strategy would result in the subdomains $X_1=[-3,-1]\times[0,1]$, $X_2=[-1,1]\times[0,1]$, $X_3=[-3,-1]\times[1,2]$ and $X_4=[-1,1]\times[1,2]$ since each dimension is partitioned into two subintervals.  
		\item \textit{Largest partitioning}: only the dimension with the largest width is partitioned into the desired number of subintervals. 
	\end{itemize}
	
	\label{para:largest_partitioning_benefit}
	Because the inclusions derived from both Natural Interval Extensions and Mean Value Form are inclusion isotonic (cf. Section \ref{subsec:intervals}), any partitioning will lead to bounds which are at least as tight as those obtained via applying interval arithmetic directly on the original domain. However, for the case of uniform partitioning we get a stronger result from Theorem~\ref{theorem:splitting}: When dividing into $N$ subintervals per dimension, the excess width associated with Mean Value Form decreases quadratically as $(w(X)/N)^2$, whereas the one given by Natural Interval Extension decreases only linearly as $w(X)/N$. Therefore, the Mean Value Form is expected to yield tighter bounds as the number of subdomains increases. 
	
	\subsubsection{The adaptive partitioning strategy}
	\label{para:partitioning_resource_dilemma}
	While the uniform partitioning provides a theoretical guarantee for bound improvement, it is not perfectly suited for fully leveraging the available GPU resources due to the lack of flexibility. Specifically, increasing the number of subintervals per dimension leads to a dramatic growth in the number of generated subdomains and as a consequence, the number of requested threads may significantly exceed the available CUDA cores on the GPU, resulting in potential thread serialization. In contrast, when decreasing the number of subintervals per dimension, much less subdomains are generated, and thus the available CUDA cores on the GPU may not be fully utilized, leading to a waste of GPU resources (see Figure~\ref{fig:CUDA_programming}). 
	
	For instance, consider partitioning a five-dimensional function on a GPU with 2560 CUDA cores available. The uniform partitioning will partition each dimension into four sub-intervals, resulting in 1024 subdomains ($4^5=1024$). Further partitioning each dimension uniformly might not be desirable because the number of created subdomains ($5^5=3125$) will be much larger than the number of available CUDA cores, potentially leading to serialization. However, if only 1024 threads are requested via CUDA, the 2560 GPU cores will not be fully used — more than half of them will be left idle. 
	
	To address the limitation of uniform partitioning, we propose an \textit{adaptive partitioning} strategy that combines the two basic partitioning strategies introduced above. In adaptive partitioning, the number of subintervals on each dimension is determined based on its relative width. Specifically, it starts from the uniform partitioning with the largest possible number of subintervals per dimension that does not yet surpass the number of available CUDA cores, and then increments the number of subintervals by one for the dimensions with relatively larger width, until the number of total subdomains is close to the number of CUDA cores available on the GPU. For example, in the case of a five-dimensional function, using adaptive partitioning, the four dimensions of largest width will be divided into five subintervals and the remaining one dimension is divided into four subintervals, which yields $5^4\cdot4 = 2500$ subdomains. Although this configuration does not fully saturate the 2560 available CUDA cores either, it significantly reduces idle resources compared to uniform partitioning, thereby improving computational efficiency. 
	
	\subsection{Implementation}
	We integrate the proposed GPU-parallel Subdomain Lower Bounding method into our open-source deterministic global optimization solver MAiNGO \cite{S1_bongartz_2018_maingo}. It is built around the B\&B algorithm, and by default uses McCormick relaxations \cite{mccormick1976computability, mitsos2009mccormick, tsoukalas2014multivariate} for bounding operations. In MAiNGO, the Directed Acyclic Graph (DAG) \cite{schichl2005interval} data structure of MC++ \cite{chachuat2015set} is leveraged to represent optimization problems. The objective function and constraints of an optimization problem are represented by different sub-graphs of the same DAG, thus only one graph is constructed at the beginning of the B\&B algorithm. 
	
	By default, interval extensions in MAiNGO (either as basis for McCormick relaxations, or if explicitly chosen as bounding method) are computed on the CPU, by evaluating the DAG of MC++ with the interval arithmetic of the filib++ library \cite{lerch2006filib++}. However, since these were not designed for GPU computations, dedicated implementation is needed. Specifically, in C++, we first define an \textit{interval structure} which contains two members — the lower and upper bounds of an interval, then implement the inclusions of intrinsic functions by overloading the operators and basic mathematical functions for this newly defined interval type.   
	
	However, special care is needed to ensure correct rounding behavior. While in CUDA, basic operations like $+,-,*,/,\sqrt{\cdot},$ and Fused Multiply-Add ($\texttt{fma}$) support all four rounding modes (round to nearest, round towards zero, round towards positive-infinity, and round towards negative-infinity) intrinsically, many other intrinsic functions such as $\text{pow}(\cdot)$, $\cos(\cdot)$, and $\exp(\cdot)$ do not fully comply with the accuracy requirements specified by the IEEE-754 floating-point standard. In the worst case, they may exhibit errors of up to 1–3 ulp (Unit in the Last Place) with no information indicating in which direction these errors occur. Therefore, implementations of intrinsic functions that are not internally provided by CUDA must perform double-direction rounding such that the correctness can be ensured. To address this rounding issue involved in interval evaluations, we develop the CuInterval library \cite{Kichler2025cuinterval}, offering CUDA kernels that guarantee the inclusion of the true interval result. Although double-rounding may lead to slightly wider intervals, the operations and intrinsic functions in CuInterval are carefully implemented to minimize over-approximation. A full list specifying the supported operations and their worst case over-approximations is provided in the documentation of CuInterval.
	
	In addition, since derivative information of the objective function and constraints is required by the Mean Value Form (cf. Definition~\ref{def:MVF}), we leverage  \textit{Tangent-mode (also known as Forward-mode) Algorithmic Differentiation} \cite{griewank2008evaluating} to calculate the derivatives. Analogous to the way we implement interval arithmetic, here we again define a \textit{tangent structure} which contains two members — the function value and its derivative, then overload the operators and basic mathematical functions for this newly defined tangent type. To enable derivative computations on GPUs, we develop the CuTangent library \cite{Kichler2025cutangent}, keeping seamless interoperability with the CuInterval library in mind.   
	
	\label{para:two_implementations}
	With CUDA, we develop Steps 1 and 2 of Subdomain Lower Bounding on the GPU in two implementations: the \emph{Single Kernel} implementation, and the \emph{CUDA Graph} implementation. Both implementations have two C++ classes to manage configuration parameters and memory allocation. These classes ensure that memory is correctly allocated and accessible to GPU threads across varying configurations of the Subdomain Lower Bounding method, such as the number of subdomains and the method used for constructing inclusion functions, etc.. When Mean Value Form is employed, additional memory is managed to store the derivative information required for its evaluation. Other key building blocks, including domain partitioning strategies and the evaluation of the Mean Value Form, are implemented as device functions which can be called directly on the GPU.
	
\begin{figure}[htbp]
		\centering
		\includegraphics[width=0.45\linewidth,trim = 10cm 7cm 10cm 7cm, clip]{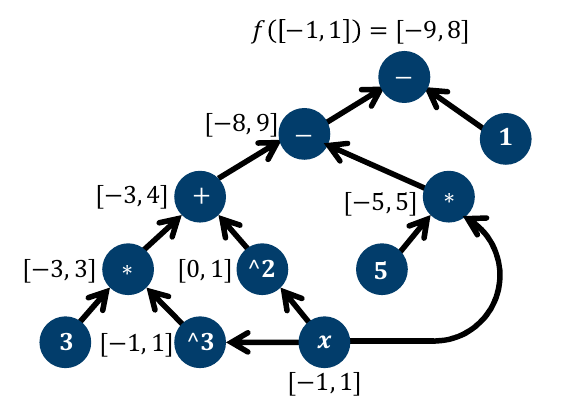}
		\caption{\footnotesize A Directed Acyclic Graph can be used to represent factorable functions: The constants and dependent variables are denoted by the leave nodes, while the intrinsic functions are represented as the inner nodes. All the nodes are connected via dependencies. Once the graph is specified, it can be evaluated in any arithmetic, including interval arithmetic to get interval bounds of the corresponding function(s), here, e.g., applying Natural Interval Extension to the function $f(x)=3x^3+x^2-5x-1$.}
		\label{fig:DAG_NIE}
	\end{figure}

	\label{para:single_kernel}
	In the Single Kernel implementation, both the domain partitioning and parallel lower bounding are wrapped into the same CUDA kernel, which will then be executed simultaneously by all the requested GPU threads if the GPU resources allow. During the parallel lower bounding, a simpler DAG datastructure constructed from the MC++ DAG is evaluated over each subdomain, via either Natural Interval Extension or Mean Value Form (see Figure~\ref{fig:DAG_NIE} for the example of evaluating $f(x)=3x^3+x^2-5x-1$ via Natural Interval Extension). The data communications, i.e., memory copying operations from the host to the GPU and vice versa, are carried out via CUDA API calls before the kernel is launched and after its execution is finished, respectively.  
	
	In the CUDA Graph implementation, the domain partitioning logic is encapsulated within a dedicated kernel, which serves as the root node of the CUDA Graph. The DAG representing the objective function and constraints is then mapped onto the same CUDA Graph structure. Each node in the DAG, corresponding to an intrinsic function, is replaced by a CUDA kernel that performs the equivalent computation. These kernels are interconnected within the CUDA Graph according to the dependency structure defined in the original DAG, such that the computation flow can be preserved. In addition, CUDA API calls are further appended to the CUDA Graph to facilitate data transfer from GPU to CPU. These calls are configured to depend on the completion of the kernels associated with the objective function and constraints, ensuring that the interval bounds are copied back to the CPU only after the relevant computations have finished. The CUDA Graph is instantiated once but launched every time the bounding operation is performed. During each launch, the kernels on the CUDA Graph are scheduled to all the requested GPU threads and thus are executed in parallel.  
	
	\label{para:CPU_serial_version}
	To enable analyzing the effect of the partitioning into subdomains in isolation from the GPU parallelization, we also implement a serial version of Subdomain Lower Bounding, where the subdomains are evaluated in serial, i.e., within a for-loop on the CPU side. By setting the number of subdomains as one, it reduces to a single-interval benchmark, where interval arithmetic is directly applied on the entire B\&B domain without splitting.
	
	Besides basic B\&B, MAiNGO also provides common bound-tightening techniques; however, among these techniques, currently only Constraint Propagation \cite{schichl2005interval} is available for the newly developed Subdomain Lower Bounding method, since the remaining ones are not readily compatible with interval methods but require non-constant (convex) relaxations.
	
	\section{Numerical Experiments}
	The proposed GPU-parallel Subdomain Lower Bounding, together with the other building blocks of a spatial B\&B algorithm, constitute a heterogeneous spatial B\&B framework. In this section, we  analyze the performance of this heterogeneous spatial B\&B algorithm with multiple numerical experiments. First, we focus on the benefits of using more subdomains as well as the speedups that can be achieved via GPU parallelization. Then we discuss the performance of the Single Kernel implementation and that of the CUDA Graph implementation. Finally, we compare the proposed approach with MAiNGO's default solver, which employs McCormick relaxations for bounding. 
	
	All numerical experiments are carried out in the scope of double precision (FP64), on an HP personal laptop with a 13th Gen Intel(R) Core(TM) i7-13700H CPU and an Nvidia RTX A1000 6GB Laptop GPU (2560 CUDA cores available). The Nvidia RTX A1000 Laptop GPU exhibits a substantial gap between single precision (FP32) and double precision computational throughput, with an approximate FP64/FP32 performance ratio of 1:64. This limited availability of FP64-capable hardware resources can lead to thread-level scheduling bottlenecks, i.e., GPU threads may stall while awaiting access to FP64 execution units, which will reduce the overall occupancy and thus the parallel efficiency. Despite such architectural limitations, our heterogeneous spatial B\&B algorithm still achieves notable performance improvements, as evidenced by the results presented below.
	
	\subsection{Benefits of domain partitioning for lower bounding}
	To evaluate the performance of the proposed method, we use example problems containing artificial neural networks (ANNs). ANNs can model complex nonlinear relationships and are widely used as surrogate models nowadays. However, solving problems involving trained ANN models to global optimality remains challenging.
	
	As illustrative example, we train an ANN surrogate model for the Peaks function \cite{Schweidtmann2018}
	\begin{equation*}
		f_\textup{Peaks}\left(x_1,x_2\right) := 3\left(1 - x_1\right)^2 e^{-x_1^2 - \left(x_2 + 1 \right) ^2} 
		- 10\left(\frac{x_1}{5} - x_1^3 - x_2^5\right) e^{-x_1^2 - x_2^2} - \frac{e^{-\left(x_1 + 1\right) ^2 - x_2^2}}{3},
	\end{equation*}
	on the domain $[-3, 3]^2$. The $\tanh$ function is used as a popular activation function resulting in a smooth ANN. In some sense, $\tanh$ also represents one of the worst -case scenarios for GPUs since $\tanh$ is notoriously slow on GPUs. Other network architectures could thus benefit even more from GPUs as has been observed plenty of times in regular machine learning literature (i.e., outside the context of global optimization) \cite{pmlr-v15-glorot11a}. 
	
	The trained ANN model is denoted as $\textup{ANN}_{f_\textup{Peaks}}$, and contains two hidden layers with $10$ neurons and $8$ neurons, respectively. The global minimum of the trained surrogate model can be found via
	\begin{equation}
		\label{eq:ANN_peaks}
		\min_{\textit{\textbf{x}} \in [-3, 3]^2} \textup{ANN}_{f_\textup{Peaks}}(\textit{\textbf{x}}). 
	\end{equation}  
	An even harder optimization problem can be constructed as
	\begin{equation}
		\label{eq:validate_ANN_peaks}
		\min_{\textit{\textbf{x}} \in [-3, 3]^2} \textup{ANN}_{f_\textup{Peaks}}(\textit{\textbf{x}}) - f_\textup{Peaks}(\textit{\textbf{x}}), 
	\end{equation}
	which together with the analogous maximization problem can serve to bound the error of the ANN  surrogate model.
	
\begin{figure}[htbp]
		\centering
		\begin{subfigure}[t]{0.45\textwidth}
			\centering
			\includegraphics[width=\linewidth]{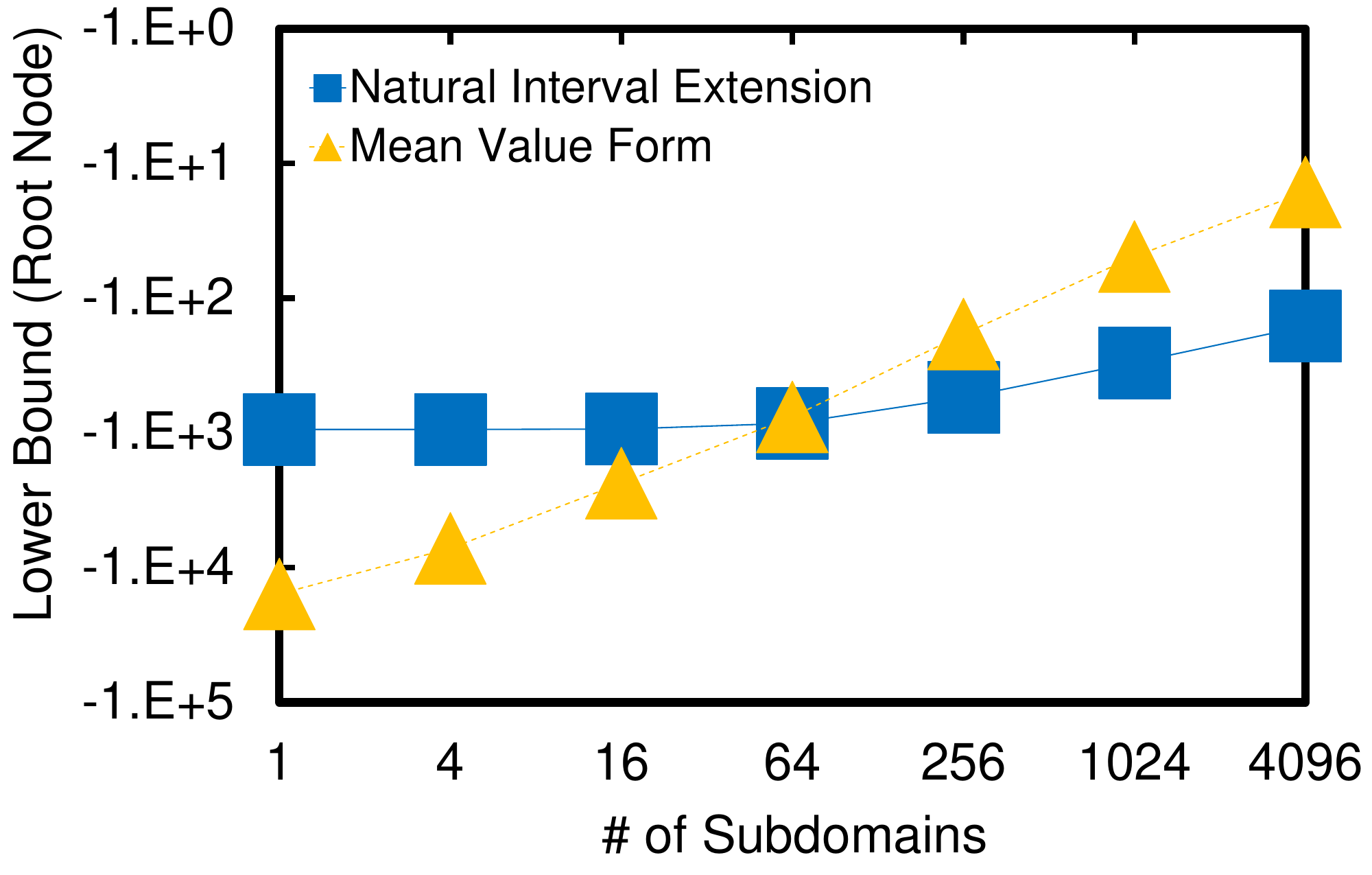}
			\vspace{-3mm}
			\caption{}
			\label{fig:LB_subdomains_peakANN}
		\end{subfigure}
		\vspace{5mm}
		\hfill
		\begin{subfigure}[t]{0.45\textwidth}
			\centering
			\includegraphics[width=\linewidth]{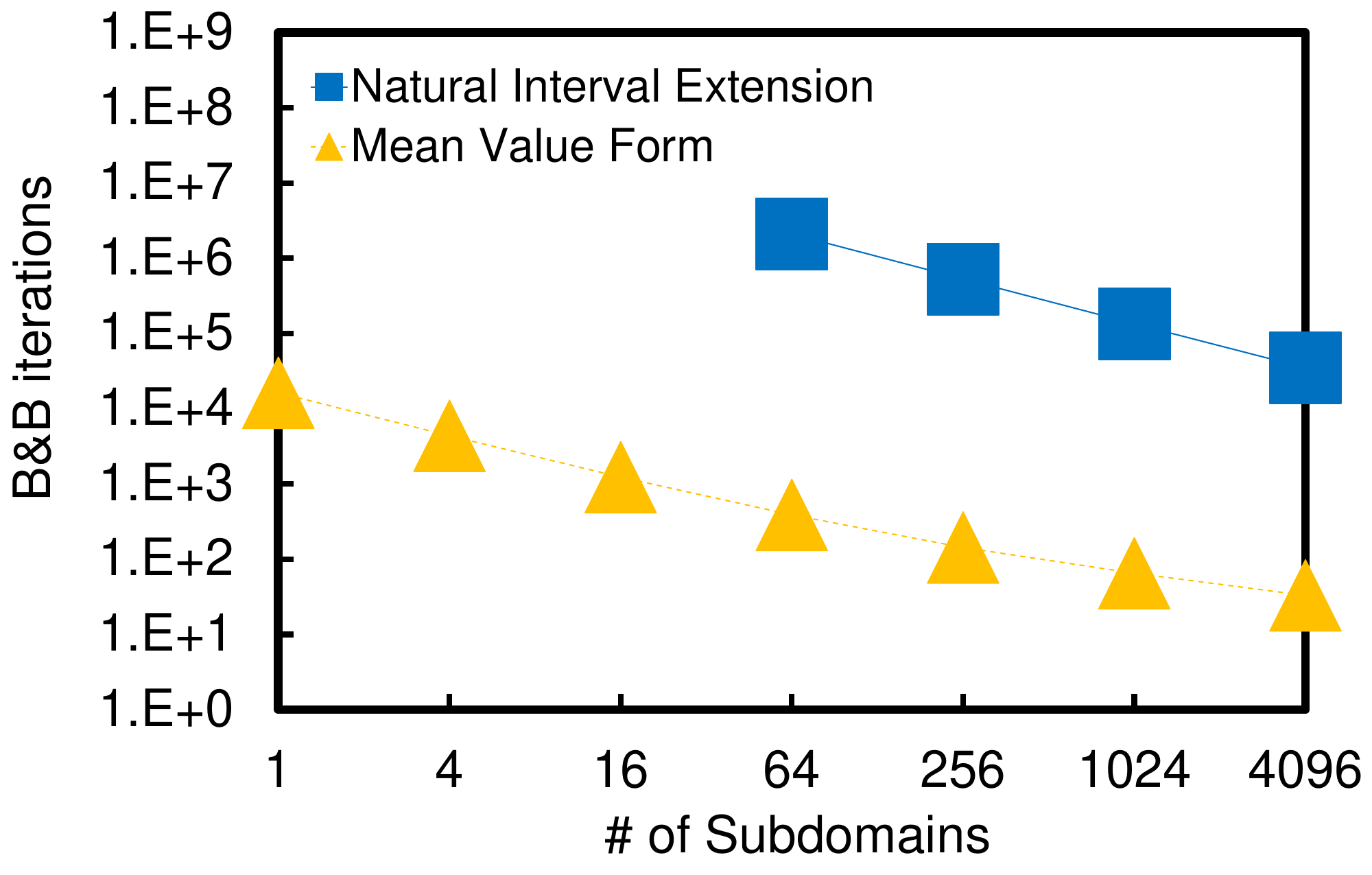}
			\vspace{-3mm}
			\caption{}
			\label{fig:BBiterations_subdomains_peakANN}
		\end{subfigure}
		\vspace{5mm}
		\begin{subfigure}[t]{0.45\textwidth}
			\centering
			\includegraphics[width=\linewidth]{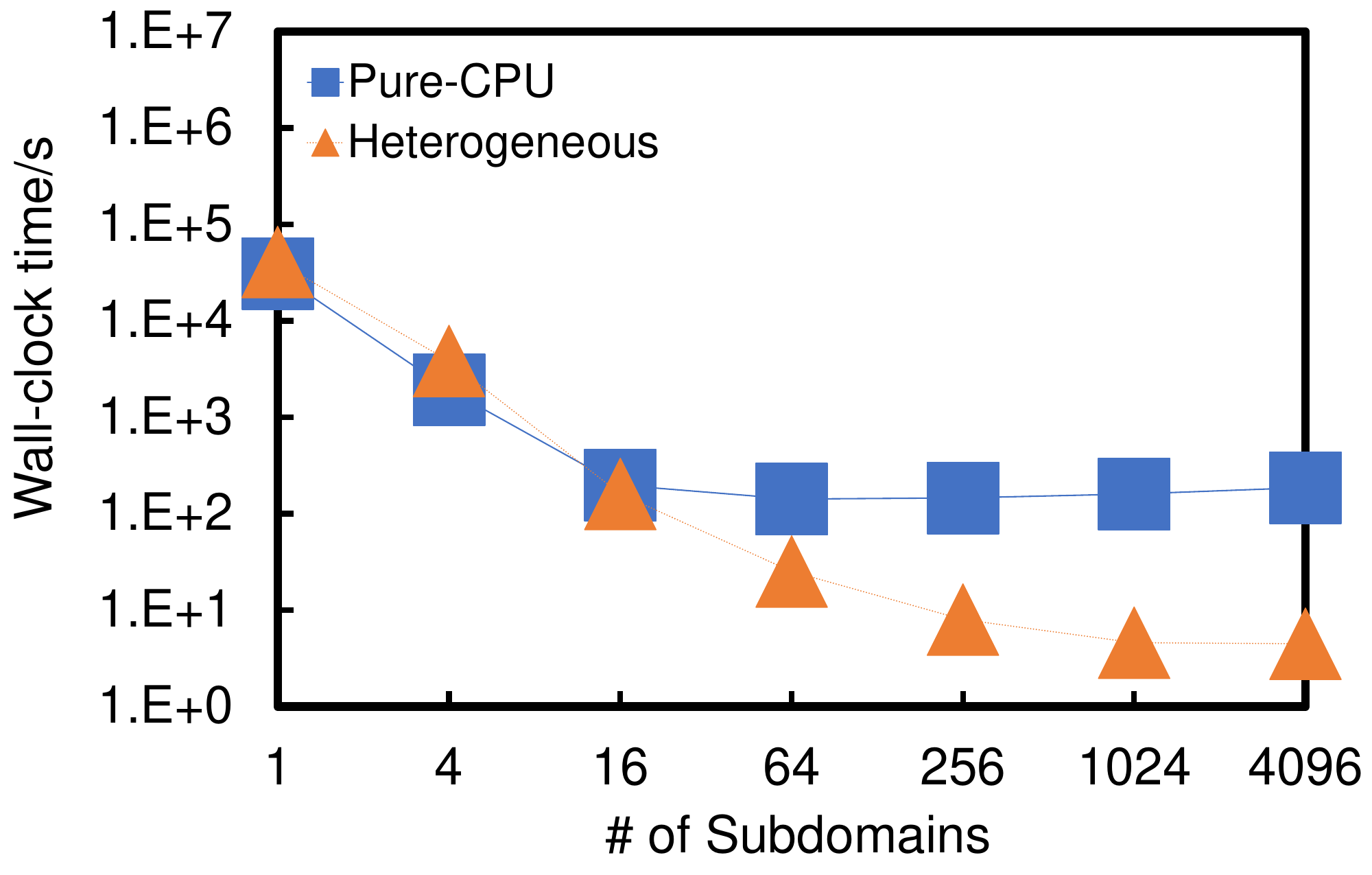}
			\vspace{-3mm}
			\caption{}
			\label{fig:Wallclock_serialParallel_modelvaliPeakANN}
		\end{subfigure}
		\hfill
		\begin{subfigure}[t]{0.45\textwidth}
			\centering
			\includegraphics[width=\linewidth]{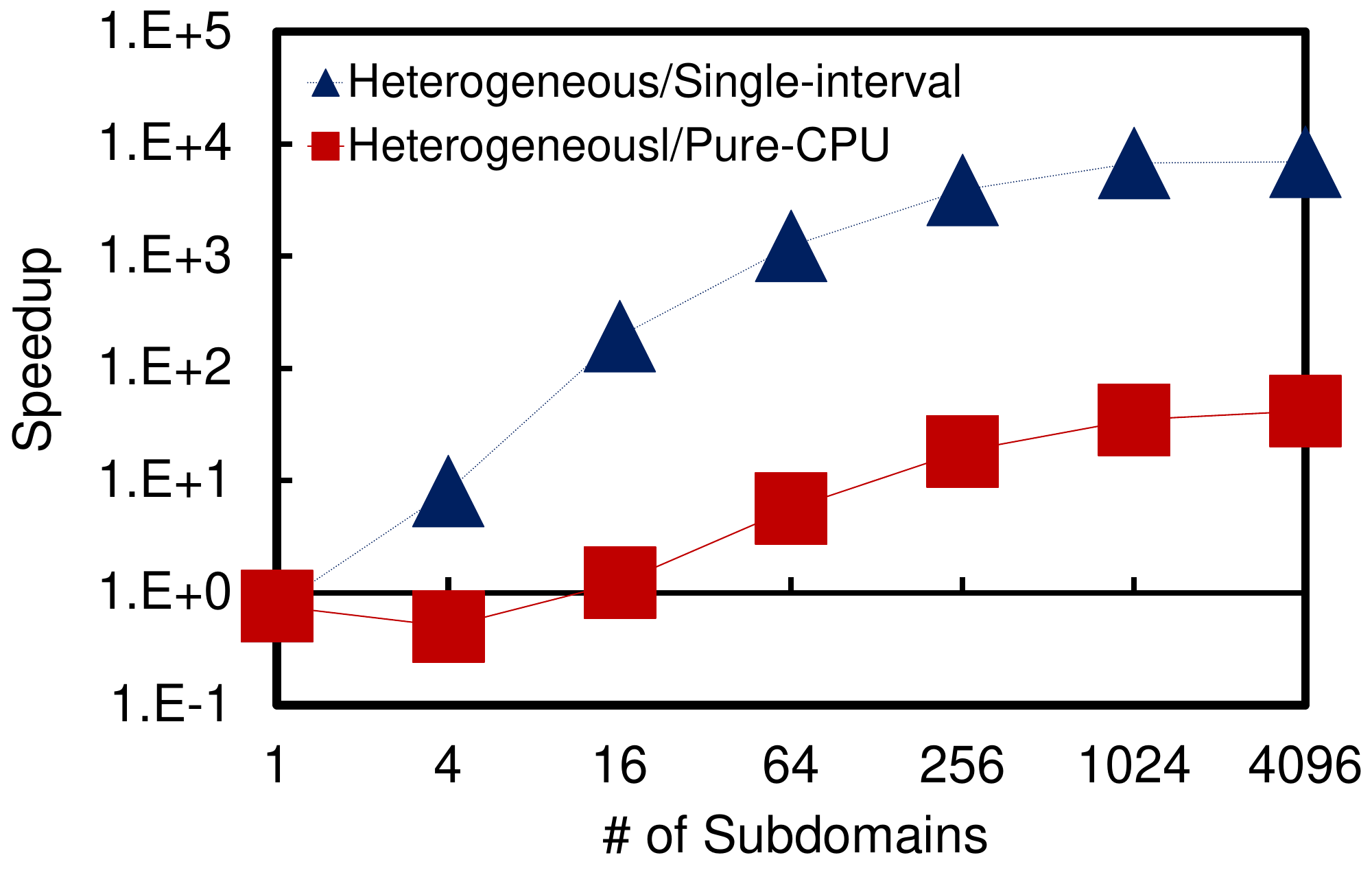}
			\vspace{-3mm}
			\caption{}
			\label{fig:Speedup_serialParallel_modelvaliPeakANN}
		\end{subfigure}
		\caption{\footnotesize Using more subdomains per B\&B node makes the lower bounds tighter, see (a) for an example of the B\&B root node. Tighter lower bounds result in less B\&B iterations, on which Mean Value Form outperforms Natural Interval Extension for all given numbers of subdomains (b). Moreover, GPU parallelization makes it possible to evaluate numerous subdomains simultaneously, therefore lower wall-clock time is achieved (c), and the spatial B\&B algorithm reaches higher speedups (d) as the number of subdomains increase. To make sure the Natural Interval Extension converges within the time limit, the tests in subfigures (a) and (b) are carried out on Problem~\eqref{eq:ANN_peaks}; In (c) and (d), only Mean Value Form is used, and the tests are carried out using the CUDA Graph implementation on Problem~\eqref{eq:validate_ANN_peaks}, which is more complex and thus visualizes scaling on a more relevant problem. In all these tests, no bound tightening technique is applied (i.e., pure B\&B).} 
		\label{fig:benefits_of_subdomains}
	\end{figure}

	We first evaluate the effect of the number of subdomains in our Subdomain Lower Bounding method on Problem~\eqref{eq:ANN_peaks}, where the subdomains are partitioned according to the uniform partitioning strategy. As anticipated, increasing the number of subdomains leads to much tighter lower bounds (see Figure~\ref{fig:LB_subdomains_peakANN} for the example of the root node). The lower bounds computed with Mean Value Form are less tight than those obtained via Natural Interval Extension when a small number of subdomains is employed. However, the Mean Value Form benefits more from the splitting technique (see Section~\ref{para:excesswidth}), such that its bounds get tighter than those from Natural Interval Extension as more subdomains are used.
	
	The improvement in bound tightness facilitates earlier pruning of non-promising nodes, thereby reducing the total number of iterations required for the convergence of B\&B (see Figure~\ref{fig:BBiterations_subdomains_peakANN}).
	In terms of B\&B iterations, the Mean Value Form always outperforms the Natural Interval Extension, even though its root node relaxation was weaker when using few sobdomains. This can be attributed to its higher convergence order \cite{hansen1969topics}, which results in sharper bounds as the width of the B\&B nodes shrinks through branching.
	
	According to previous analysis, more subdomains enable less B\&B iterations, but it is not clear whether this will reduce computational time, given that each B\&B iteration now requires bounding on more subdomains, and transferring data between devices to enable GPU parallelization. To analyze the speedups that can be reached with our heterogeneous spatial B\&B framework, we carry out the comparison against the pure-CPU B\&B algorithm in which the serial implementation of the Subdomain Lower Bounding method (as explained in Section \ref{para:CPU_serial_version}) is leveraged for lower bounding. 
	
	As illustrated in Figure~\ref{fig:Wallclock_serialParallel_modelvaliPeakANN}, we first observe that even the serial pure-CPU implementation of the Subdomain Lower Bounding procedure benefits from more subdomains: Using up to 64 subdomains reduces the wall-clock time by approximately two orders of magnitude compared to no domain partitioning (i.e., using one subdomain without splitting). This indicates that the \emph{temporary} massive subdivision we use here is more efficient than fully relying on the sequential univariate subdivision through branching in the B\&B tree up to this point. As mentioned in Section \ref{subsec:intervals}, the temporary subdivision used herein bears some resemblance to piece-wise linear relaxations that have been employed for some problem types \cite{hasan2010piecewise,misener2011apogee} (though not in the context of parallelization), in which B\&B nodes were also subdivided to obtain tighter bounds. However, while these methods solved a mixed-integer problem to compute the overall bound, here we merely have to take the interval hull, which is much faster - likely at the expense of weaker bounds.
	For our present example, the performance of the pure-CPU B\&B plateaus as the number of subdomains increases beyond 64, indicating that the reduction in B\&B iterations due to tighter bounds (achieved via Subdomain Lower Bounding) is insufficient to compensate for the overhead introduced by the additional bounding operations conducted in serial. 
	
	\label{para:performance_saturation}
	In contrast, the heterogeneous implementation continues to achieve reductions in wall-clock time by up to over three orders of magnitude until the number of subdomains reaches 1024, Although further increasing the number of subdomains to 4096 yields tighter bounds and fewer B\&B iterations (cf. Figure~\ref{fig:LB_subdomains_peakANN} and \ref{fig:BBiterations_subdomains_peakANN}), the resulting overall wall-clock time remains on a similar level to that achieved with 1024 subdomains. This saturation in performance can be attributed to two factors. First, the GPU used in our experiments is equipped with 2560 CUDA cores, while 4096 threads are requested — one per subdomain. Since the number of threads exceeds the available cores and the double precision performance is constrained by the FP64/FP32 ratio, full parallel execution cannot be realized, leading to degradation in efficiency. Second, when computing the refinement in Step 3 of the Subdomain Lower Bounding method (cf. Section~\ref{para:steps_on_GPU}), a serial sorting is performed on the CPU side. As the number of subdomains increases, more comparison steps have to be done which also contributes to longer processing times.
	
	Despite the overhead associated with memory transfers and the generally higher per-core performance of CPUs compared to individual GPU cores, the heterogeneous spatial B\&B algorithm still achieves significant speedups over its pure-CPU counterpart when using more subdomains, thanks to the parallel execution on the GPU (see Figure~\ref{fig:Wallclock_serialParallel_modelvaliPeakANN}, and red squares in Figure~\ref{fig:Speedup_serialParallel_modelvaliPeakANN}). Moreover, we further compare the heterogeneous B\&B with the single-interval benchmark (cf. Section~\ref{para:CPU_serial_version}.), i.e., applying interval arithmetic without splitting purely on the CPU, to evaluate the combined benefits originating from the tighter bounds obtained via splitting technique, and the computational efficiency offered by GPU parallelization. Remarkably, the heterogeneous approach attains speedups over three orders of magnitude against the single-interval benchmark (see blue triangles in Figure~\ref{fig:Speedup_serialParallel_modelvaliPeakANN}). The speed gain exceeds the number of subdomains employed, consistent with our expectation that combining the splitting technique with GPU parallelization enables a lower-bounding solver that yields tighter bounds and meanwhile remains efficient in evaluation.
	
	\subsection{Single Kernel vs. CUDA Graph}
	As described in Section~\ref{para:two_implementations}, the GPU part of the proposed method is developed into two implementations: the Single Kernel one and the CUDA Graph one. In the Single Kernel implementation, the DAG representing the objective function and constraints resides on the GPU (cf. Section~\ref{para:single_kernel}). Since all GPU tasks are wrapped into one kernel, the DAG nodes are evaluated sequentially, even though some of them may be independent from each other and thus suitable for parallel evaluation. In contrast, in the CUDA Graph implementation, each DAG node corresponding to an intrinsic function, is mapped as an individual kernel within the CUDA graph. This design enables the parallel execution of independent operations or functions, provided that sufficient GPU resources are available (cf. Figure~\ref{fig:CUDA_graph}).
	
\begin{figure}[htbp]
		\centering
		\begin{subfigure}[t]{0.45\textwidth}
			\centering
			\includegraphics[width=\linewidth]{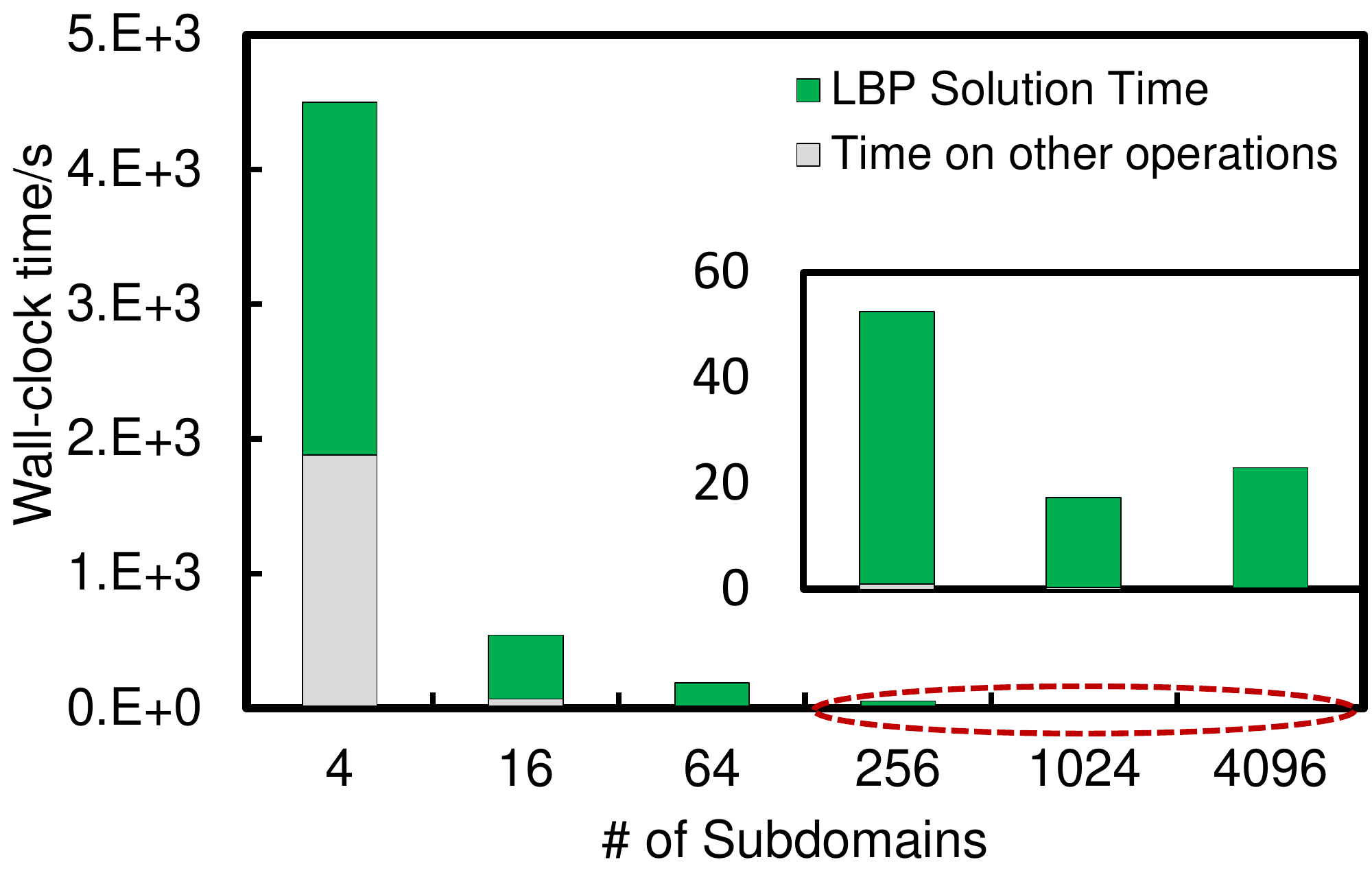}
			\vspace{-3mm}
			\caption{}
			\label{fig:timecomposition_modelvali_sk}
		\end{subfigure}
		\hfill
		\begin{subfigure}[t]{0.45\textwidth}
			\centering
			\includegraphics[width=\linewidth]{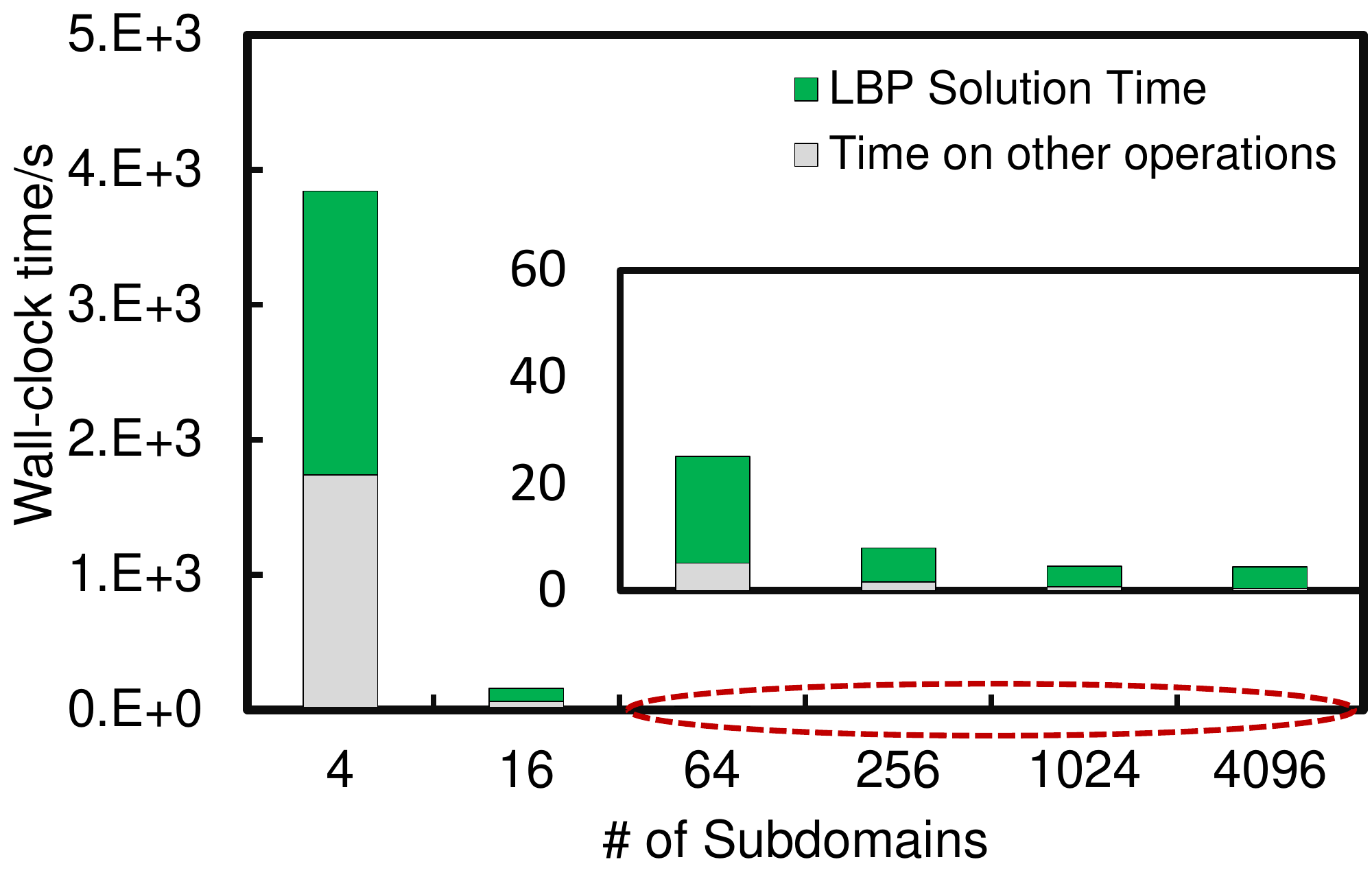}
			\vspace{-3mm}
			\caption{}
			\label{fig:timecomposition_modelvali_cu}
		\end{subfigure}
		\caption{\footnotesize Compared to the Single Kernel implementation (a), the CUDA Graph version (b) attains lower wall-clock time since it spends much less time on solving the lower bounding problem (LBP). The inset figure shows a zoomed-in version of the red-dashed ellipse. The test is carried on Problem~\eqref{eq:validate_ANN_peaks} and no bound tightening technique is applied.} 
		\label{fig:implementation_comparison}
	\end{figure}

	As shown in Figure~\ref{fig:timecomposition_modelvali_sk} and \ref{fig:timecomposition_modelvali_cu}, compared to the Single Kernel implementation, the CUDA Graph version reduces the lower bounding time by approximately 20 to 90\%, depending on the number of subdomains. This higher efficiency in executing the proposed Subdomain Lower Bounding method may be attributed to the (potential) more parallelism offered by the CUDA graph and the more efficient memory management. Note that the reduction in time can be problem-dependent: In preliminary experiments with simpler problems, we observed that the two implementations could offer similar performance. 
	
	\subsection{Comparison with McCormick-based Solver}
	Since many current global solvers use convex relaxations instead of intervals, we compare the proposed heterogeneous spatial B\&B algorithm with MAiNGO's default solver which runs in serial purely on the CPU and employs (multivariate) McCormick relaxations \cite{mccormick1976computability, tsoukalas2014multivariate, mitsos2009mccormick} for bounding. The tests are carried out on problems of different characteristics, and the effect of bound tightening techniques is also investigated.
	
	\subsubsection{Scalability}
	As discussed in Sections~\ref{para:partitioning_resource_dilemma} and \ref{para:performance_saturation}, more serialization occurs when the number of requested threads keeps growing, which imposes a practical constraint on the number of subintervals that we can partition along each dimension. This will in turn limit the sharpness of the refinement, leading to a degradation on the performance of B\&B. This limitation could be expected to be more relevant for problems of higher dimensionality. Therefore, to get some first insights into the scalability of the proposed heterogeneous spatial B\&B framework, we evaluate its performance on problems with different numbers of variables and compare it against MAiNGO's default solver. Note that no bound tightening is used, i.e., both the Subdomain Lower Bounding method and the one based on McCormick relaxations are used in a pure B\&B configuration.
	
	Consider the $Alpine02$ function
	\begin{equation*}
		f_\textup{Alpine02,d}(x):=\prod_{i=1}^{d}\sqrt{x_i}\sin(x_i).
	\end{equation*}
	The dimensionality of the $Alpine02$ function can be adjusted by changing $d$. We train ANN surrogate models for $Alpine02$ functions with dimensions varying from two to six, on the domain $[3, 9]^d$, respectively. Each ANN model has only one hidden layer with 60 neurons. As in Section 4.1, we again consider a similar problem
	\begin{equation}
		\max_{\textit{\textbf{x}} \in [3,9]^d} \textup{ANN}_{f_\textup{Alpine02,d}}(\textit{\textbf{x}}).
		\label{eq:ANN_Alpine02}
	\end{equation}
	
	\label{para:scalability_test}
	In our heterogeneous spatial B\&B framework, we first employ both the uniform and adaptive branching strategies to get a range of reasonable numbers of subdomains, and then evaluate the method with each of these configurations. The performance yielded by using different numbers of subdomains are summarized in the Appendix (see Table~\ref{tab:Alpine02ANN_Scaling}). Note that only rather small differences in wall-clock time are observed among these different configurations.
	
\begin{figure}[htbp]
		\centering
		\begin{subfigure}[t]{0.45\textwidth}
			\centering
			\includegraphics[width=\linewidth]{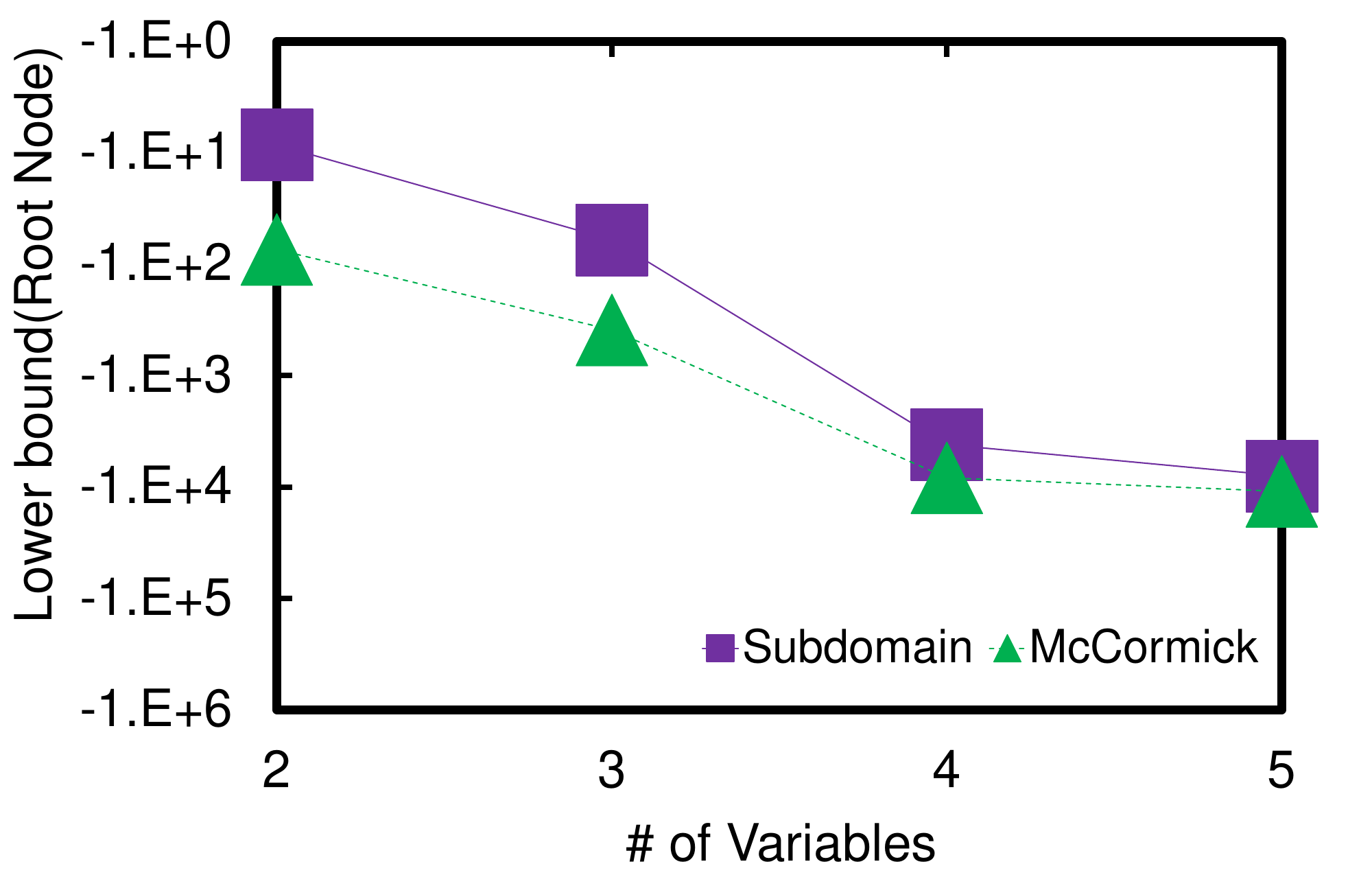}
			\vspace{-3mm}
			\caption{}
			\label{fig:Lowerbounds_alipne_scale}
		\end{subfigure}
		\hfill
		\begin{subfigure}[t]{0.45\textwidth}
			\centering
			\includegraphics[width=\linewidth]{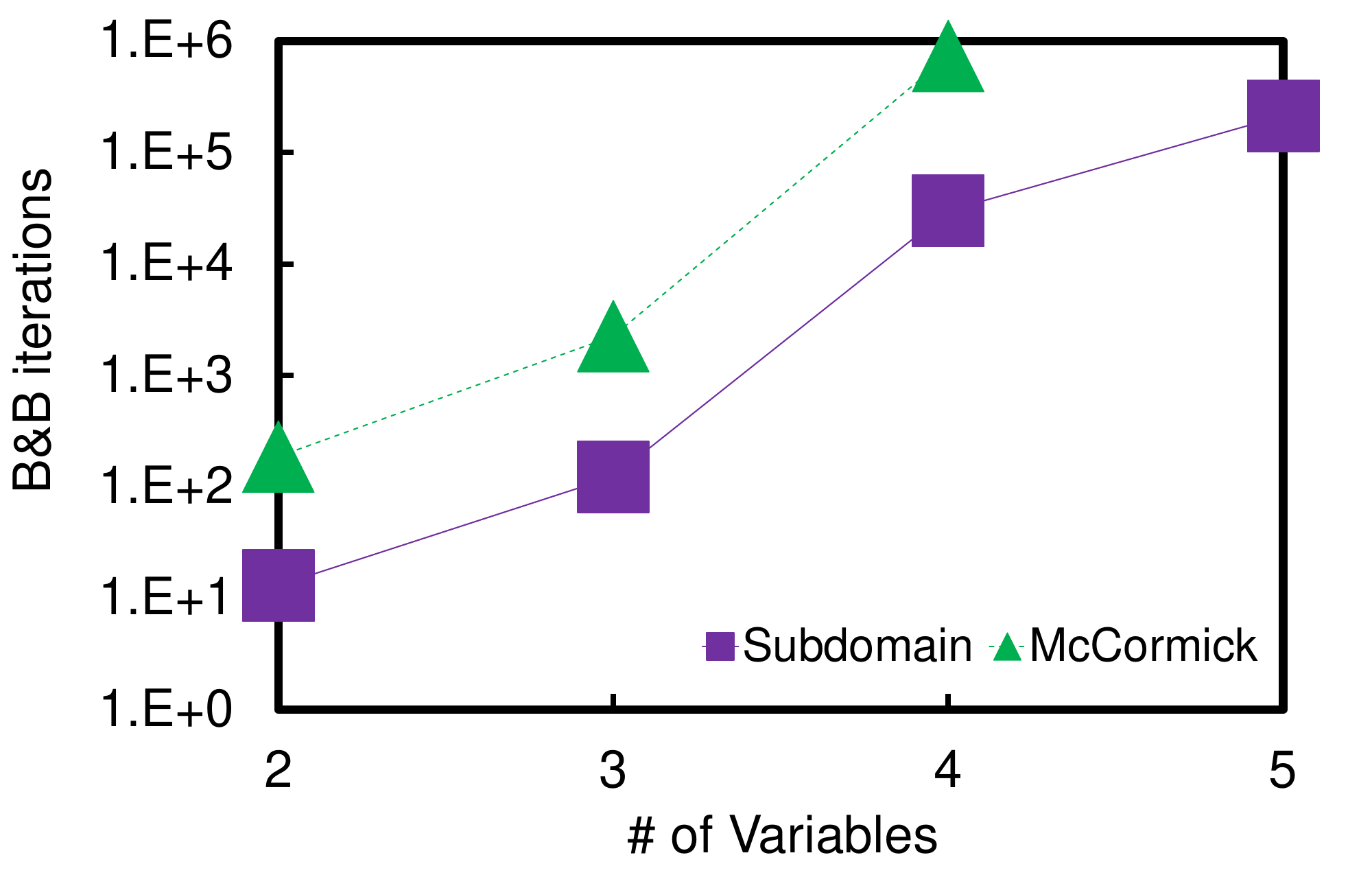}
			\vspace{-3mm}
			\caption{}
			\label{fig:BBiterations_alipne_scale}
		\end{subfigure}
		\begin{subfigure}[t]{0.45\textwidth}
			\centering
			\includegraphics[width=\linewidth]{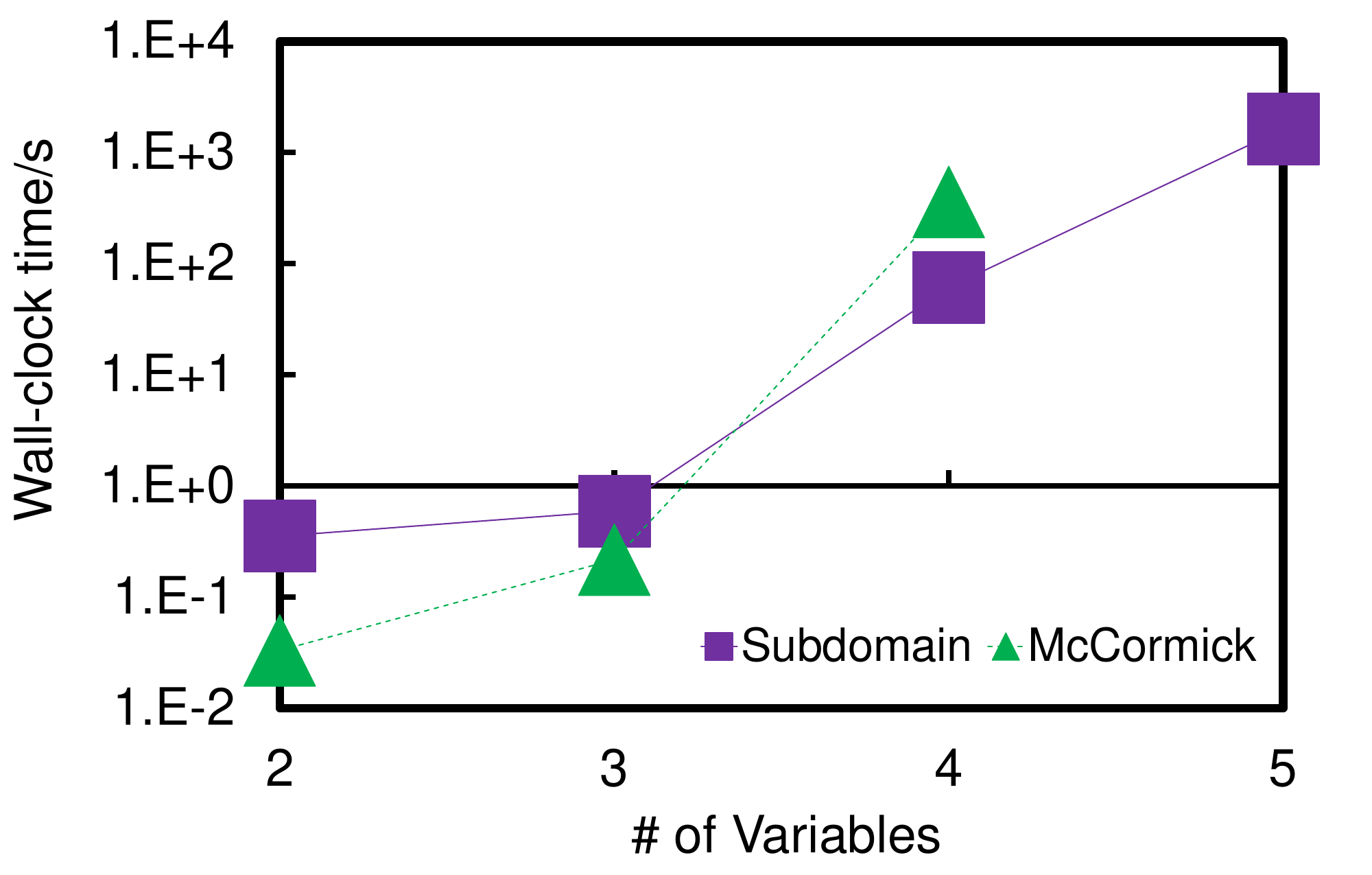}		
			\vspace{-3mm}	
			\caption{}
			\label{fig:wallclock_alipne_scale}
		\end{subfigure}
		\hfill
		\begin{subfigure}[t]{0.45\textwidth}
			\centering
			\includegraphics[width=\linewidth]{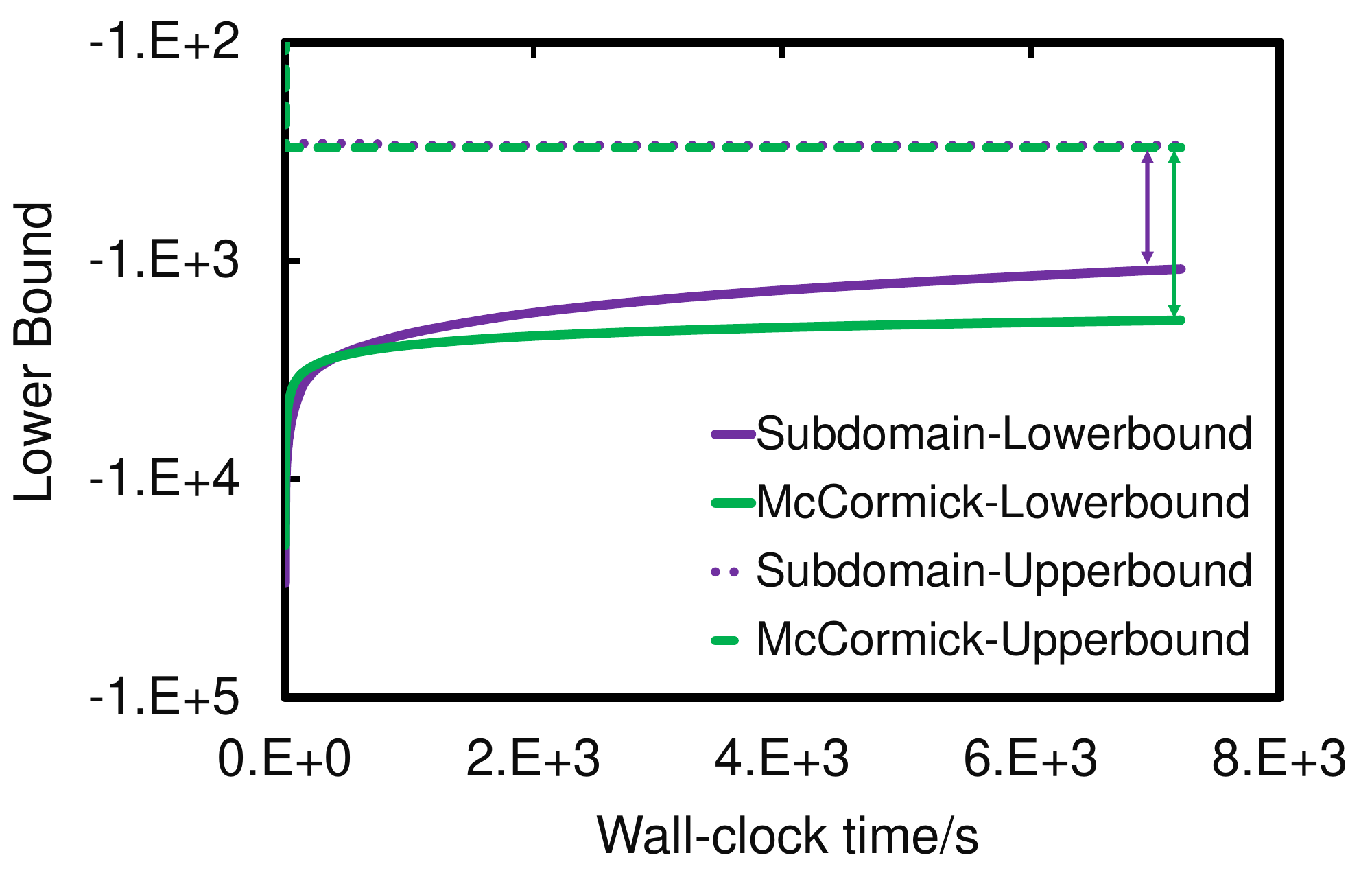}
			\vspace{-3mm}
			\caption{}
			\label{fig:Lowerbounds_alpine_6d}
		\end{subfigure}
		\caption{\footnotesize The scalability test is carried out on Problem~\eqref{eq:ANN_Alpine02} with problem dimensions varying from two to six. Here ``Subdomain" and ``McCormick" denote the results of our heterogeneous spatial B\&B leveraging the GPU-parallel Subdomain Lower Bounding method, and MAiNGO's default solver using the McCormick relaxations for bounding on the CPU, respectively. Our heterogeneous spatial B\&B outperforms MAiNGO's default solver in terms of both lower bounds (a) and B\&B iterations (b) when the problem dimension is less than six. As the dimensionality of the problem increases, both methods require substantially longer time to converge (c). Though both methods cannot converge within two hours in the six-d scenario, our heterogeneous spatial B\&B algorithm still achieves smaller absolute gap (d). In all the tests, no bound-tightening technique is applied. For the five-dimensional case, the default solver of MAiNGO fails to prove optimality within two hours, and thus its corresponding wall-clock time and number of B\&B iterations are not reported.} 
		\label{fig:method_scalability}
	\end{figure}

	In the two- to five-dimensional cases, the proposed approach outperforms the default solver of MAiNGO on both bound tightness and B\&B iterations, as illustrated in Figures~\ref{fig:Lowerbounds_alipne_scale}–\ref{fig:BBiterations_alipne_scale}. In terms of wall-clock time (see Figure~\ref{fig:wallclock_alipne_scale}), MAiNGO's default solver demonstrates relatively better efficiency on solving Problem~\eqref{eq:ANN_Alpine02} when the number of variables is fewer than four. However, in the five-dimensional scenario, it fails to prove global optimality within the two-hour time limit. In contrast, the proposed method successfully converges within one hour. Furthermore, in the six-dimensional scenario, although neither method achieves convergence within the time limit, the heterogeneous approach still yields a substantially tighter lower bound (see Figure~\ref{fig:Lowerbounds_alpine_6d}). These results demonstrate that, as problem dimensionality increases, the heterogeneous spatial B\&B algorithm remains effective and can offer superior performance compared to MAiNGO's default solver.
	
	\subsubsection{Constraint handling}
	While the previous experiments focused on unconstrained problems, we now present a preliminary proof of concept for constrained optimization by introducing various types of constraints (linear vs. nonlinear, equality vs. inequality) to the problem
	\begin{equation}
		\min_{\textit{\textbf{x}} \in [3,9]^d} f_\textup{Alpine02,d}(\textit{\textbf{x}}) - \textup{ANN}_{f_\textup{Alpine02,d}}(\textit{\textbf{x}}).
		\label{eq:validate_ANN_Alpine02}
	\end{equation}
	in the two-dimensional setting, i.e., $d=2$ (see Table~\ref{tab:SIA_MAiNGO_constrained}). 
	
	In the context of constrained problems, interval bounds of the constraints, computed during the Subdomain Lower Bounding phase, are leveraged by the proposed method to prune non-promising B\&B nodes, see Section~\ref{para:constraint_bound_pruning}. The proposed method successfully handles all reported constraint types. Though the nature of constraints has a clear impact on the performance of lower bounding methods, the heterogeneous approach consistently outperforms MAiNGO's default solver (without bound tightening) in terms of B\&B iterations, and in some cases also in terms of wall-clock time, as summarized in Table~\ref{tab:SIA_MAiNGO_constrained}.
	\begin{table}[htbp]
		\captionsetup{position=above}
		\caption{\footnotesize Evaluation of our heterogeneous spatial B\&B on Problem~\eqref{eq:validate_ANN_Alpine02} with $d = 2$ and different types of constraints, as well as the comparison against MAiNGO's default solver that leverages McCormick relaxations for bounding on the CPU. ``Subdomain" and ``McCormick" denote the results of our method and MAiNGO's default solver, respectively. 1024 subdomains are used in the test and no bound tightening technique is applied.} 
		\centering
		\includegraphics[width=\linewidth, trim=0cm 3cm 0cm 3cm, clip]{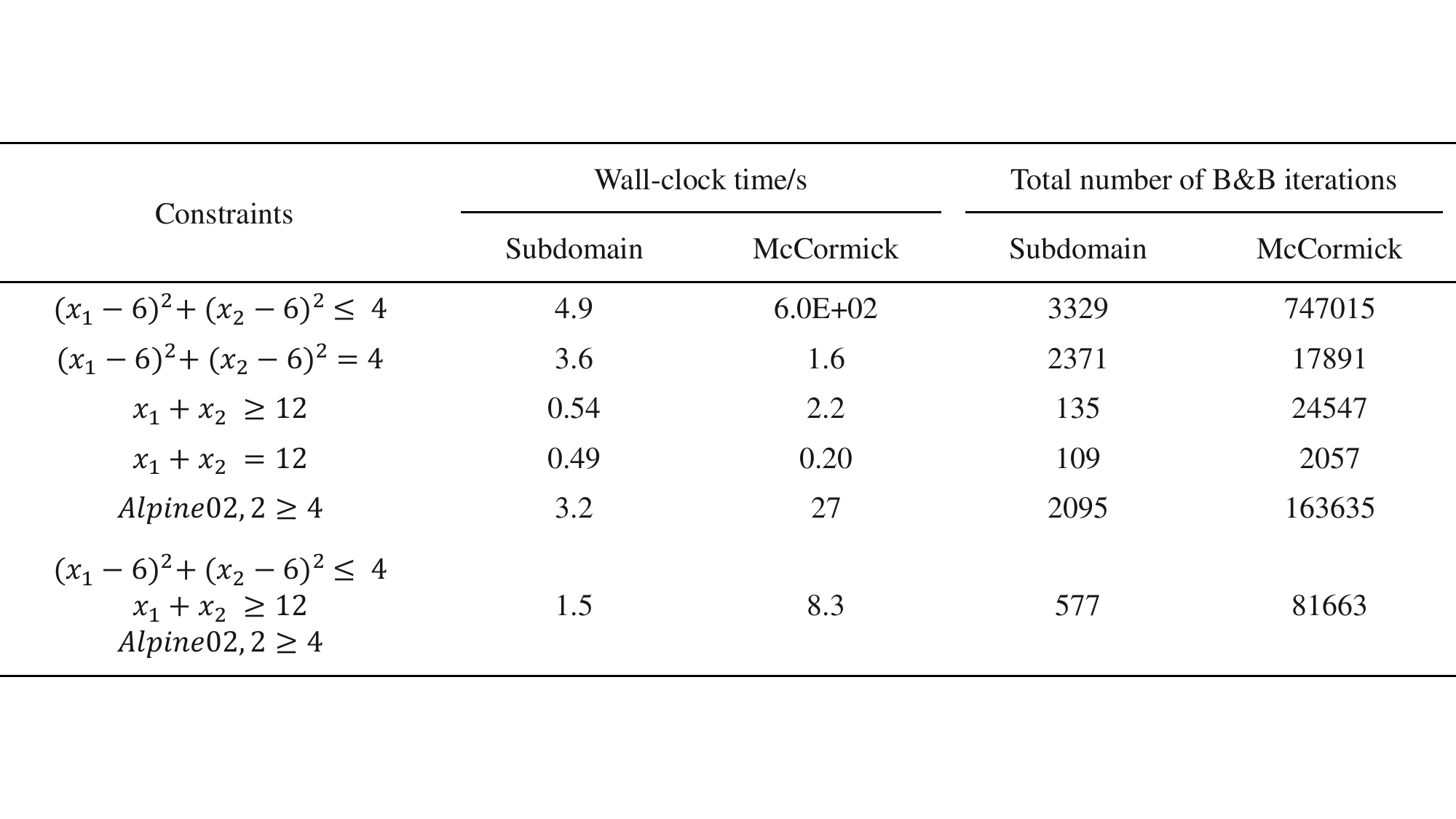}
		\label{tab:SIA_MAiNGO_constrained}
	\end{table}
	
	\subsubsection{Effect of bound tightening techniques}
	\label{para:bound_tightening_test}
	In all previous experiments, no bound tightening techniques were applied (i.e., pure B\&B). However, such techniques are commonly employed in global optimization algorithms to improve performance. Therefore, In this section, we compare the proposed heterogeneous spatial B\&B against MAiNGO's default solver, under scenarios where various bound-tightening strategies are available. Similar to Section~\ref{para:scalability_test}, the heterogeneous B\&B algorithm is evaluated with several reasonable numbers of subdomains. For each test problem used here, its definition as well as the performance yielded by using different numbers of subdomains are summarized in Appendix (see Table~\ref{tab:ComparisonProblems_details}). 
	
	Without bound tightening, the proposed heterogeneous method demonstrates competitive performance (e.g., on the Kinetic ODE Parameter Estimation problem \cite{mitsos2009mccormick}), and even significantly outperforms the default solver of MAiNGO in certain cases (e.g., on Minimizing the Styblinski–Tang function) as shown in Table~\ref{tab:SIA_MAiNGO_comparison}. Nevertheless, the default solver of MAiNGO benefits substantially from bound tightening techniques such as Duality-Based Bound Tightening (DBBT) \cite{ryoo1995global} and Constraint Propagation \cite{schichl2005interval}, which can lead to considerable speedups. Although Constraint Propagation is also applicable for the proposed method, by itself its improvement on performance is not always significant, as shown in Table~\ref{tab:SIA_MAiNGO_comparison}.
	\begin{table}[htbp]
		\captionsetup{position=above}
		\caption{\footnotesize Wall-clock time in seconds for solving different problems using either our heterogeneous spatial B\&B or MAiNGO's default solver that leverages McCormick relaxations for bounding on the CPU. Here "Subdomain" and "McCormick" denote the results of our method and MAiNGO's default solver, respectively. The last column shows the best performance that MAiNGO's default solver attains by enabling different combinations of bound-tightening techniques. For details on the Kinetic ODE parameter estimation problem, refer to \cite{mitsos2009mccormick}.}
		\label{tab:SIA_MAiNGO_comparison}
		\centering
		\includegraphics[width=\linewidth, trim=0cm 4cm 0cm 4cm, clip]{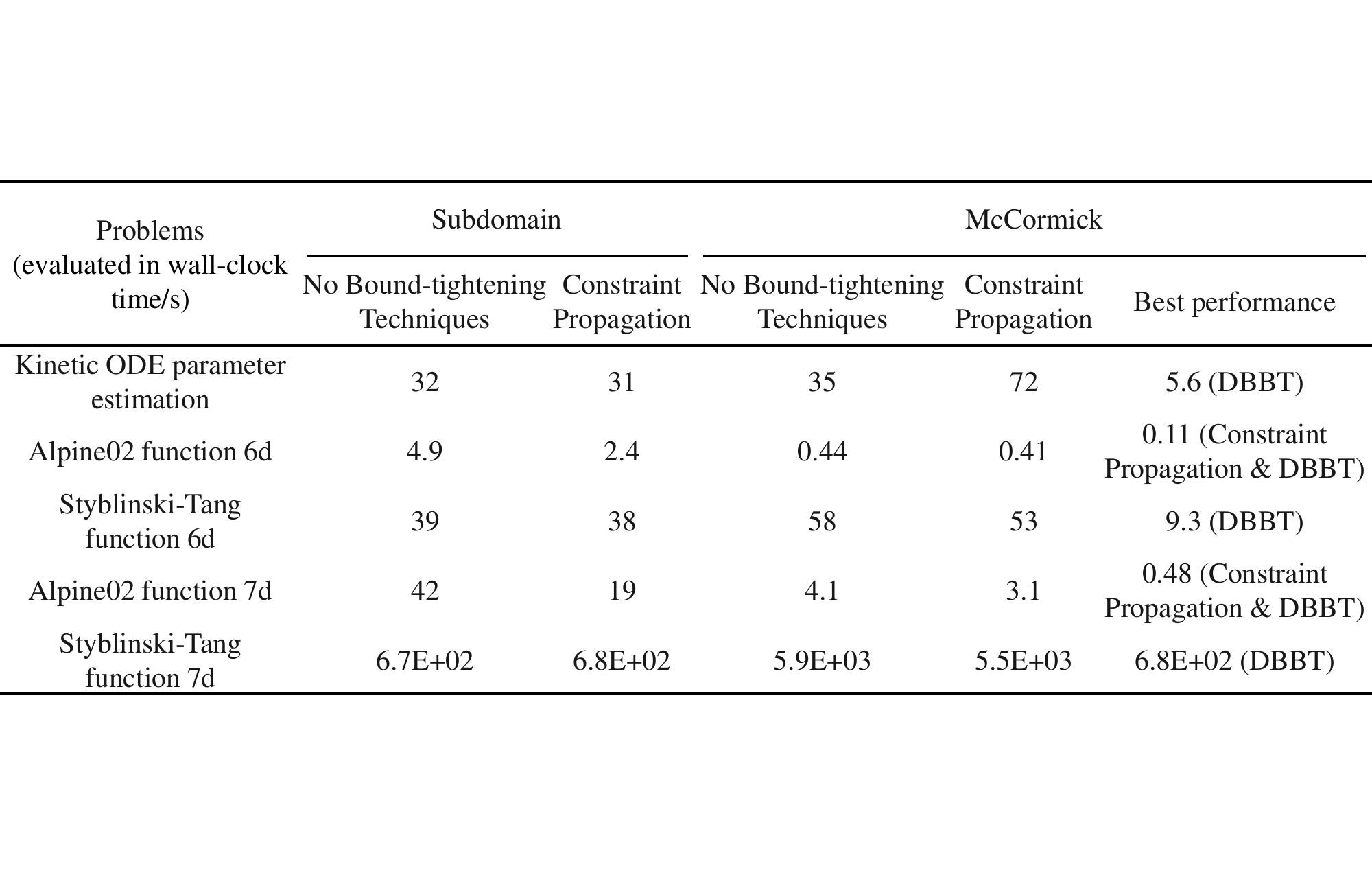}
	\end{table}
	
	\section{Conclusion and future work}
	We presented a GPU-parallel lower bounding method, termed Subdomain Lower Bounding, which can efficiently compute tight lower bounds by incorporating the splitting technique and GPU parallelization, to accelerate the spatial B\&B algorithm. 
	
	As expected, Mean Value Form yields much sharper refinements than Natural Interval Extension, benefiting from its higher convergence order. Increasing the number of subdomains leads to tighter lower bounds at B\&B nodes and consequently reduces the number of B\&B iterations. The resulting heterogeneous spatial B\&B algorithm achieves a speedup of over three orders of magnitude compared to using interval arithmetic on the CPU without splitting. Compared to the Single Kernel implementation, the CUDA Graph version demonstrates higher efficiency in solving lower bounding problems, benefiting from the better memory management and more parallelism offered by the CUDA Graph. Regarding the comparison with MAiNGO's default solver, the heterogeneous spatial B\&B delivers competitive or even better performance on some test problems, but can also be less competitive, especially when bound tightening is applied. 
	
	The results presented herein are an indication of the performance that can be achieved on a relatively accessible computing infrastructure. The proposed heterogeneous spatial B\&B algorithm is expected to deliver improved performance on hardware with more powerful GPUs. A key source of this improvement lies in the ability to process a larger number of subdomains in parallel, since high-end GPUs generally feature a significantly greater number of CUDA cores. Moreover, the GPUs engineered for HPC environments typically offer a much higher FP64/FP32 performance ratio compared to consumer-grade GPUs. For example, the NVIDIA H100 GPU has a ratio of 1:2, representing a substantial improvement over the 1:64 ratio observed in the NVIDIA RTX A1000 Laptop GPU. This enhanced support for double precision arithmetic markedly alleviates the scheduling bottlenecks associated with FP64 operations, thereby improving the overall parallel efficiency. 
	
	The gap between FP64 and FP32 performance introduces complexity in selecting an appropriate number of subdomains, because not all GPU-executed instructions perform double-precision arithmetic. In addition, the serial comparisons on the CPU also poses a penalty on using more subdomains. However, although some test cases indicate that fewer subdomains may lead to better performance, the observed differences across varying numbers of subdomains are modest for smaller-scale problems. Moreover, increasing the number of subdomains leads to improved performance for relatively larger problems (see Tables~\ref{tab:Alpine02ANN_Scaling} and \ref{tab:ComparisonProblems_details}). Thus, in general, saturating the available CUDA cores on the GPU could be a favorable strategy.
	
	While the present results demonstrate promising potential for the proposed method, several improvements could still be explored. At present we leverage the GPU exclusively for parallel lower bounding to enable the efficient evaluation of tight bounds, based on the consideration that other components of the B\&B algorithm can be better suited for CPU parallelization. Therefore, one promising direction is parallelizing the B\&B algorithm at the node level on the CPU (cf. Section~\ref{para:parallelization_types}), wherein the GPU is leveraged to perform parallel lower bounding for each B\&B node. This CPU-GPU parallelization design is expected to exploit both the CPU and GPU hardware resources especially in HPC clusters equipped with massive CPU cores and multiple powerful GPU cards. 
	
	Besides, in the current implementation, the CPU remains idle while the GPU evaluates interval bounds over subdomains, and conversely, the GPU is idle when the CPU processes other components of the B\&B algorithm. This sequential utilization of hardware resources results in a waste of the computation capacity. Exploring strategies to enable concurrent CPU-GPU execution such that the hardware utilization can be maximized, presents another interesting direction for future research. Moreover, developing bound-tightening strategies and heuristics specifically tailored to Subdomain Lower Bounding could further improve its efficiency. Finally, it is worth further evaluating the proposed heterogeneous spatial B\&B algorithm on established benchmark libraries and across computing platforms with diverse hardware configurations.
	
	\paragraph{Acknowledgement} 
	Hongzhen Zhang and Dominik Bongartz gratefully acknowledge funding through Internal Funds KU Leuven (STG/22/060). Tim Kerkenhoff, Manuel Dahmen and Alexander Mitsos received funding from the Helmholtz Association of German Research Centers. The authors thank Clara Witte for her support in integrating the code into MAiNGO.
	
	\printbibliography	
	
	\newpage
	\appendix
	\section{Appendix}
	We summarize the optimization problems leveraged to test the proposed method. As discussed in Section~\ref{para:performance_saturation}, more subdomains may not always lead to less computation time due to the existence of serialization. Since the resulting performance can be problem-dependent, for each optimization problem, we applied multiple reasonable numbers of subdomains that are generated via either \textit{uniform partitioning} or \textit{adaptive partitioning}. 
	\begin{table}[htbp]
		\captionsetup{position=above}
		\caption{\footnotesize Details on the scalability test of the proposed heterogeneous spatial B\&B algorithm (used in Section~\ref{para:scalability_test}). This test is carried out on Problem~\eqref{eq:ANN_Alpine02} with two to six variables. The last column shows the performances (regarding wall-clock time) attained at different numbers of subdomains. For the six-dimensional scenario, the final absolute gap is reported since the method fails to converge within two hours.}
		\centering
		\includegraphics[width=\linewidth, trim=0cm 0.5cm 0cm 0.5cm, clip]{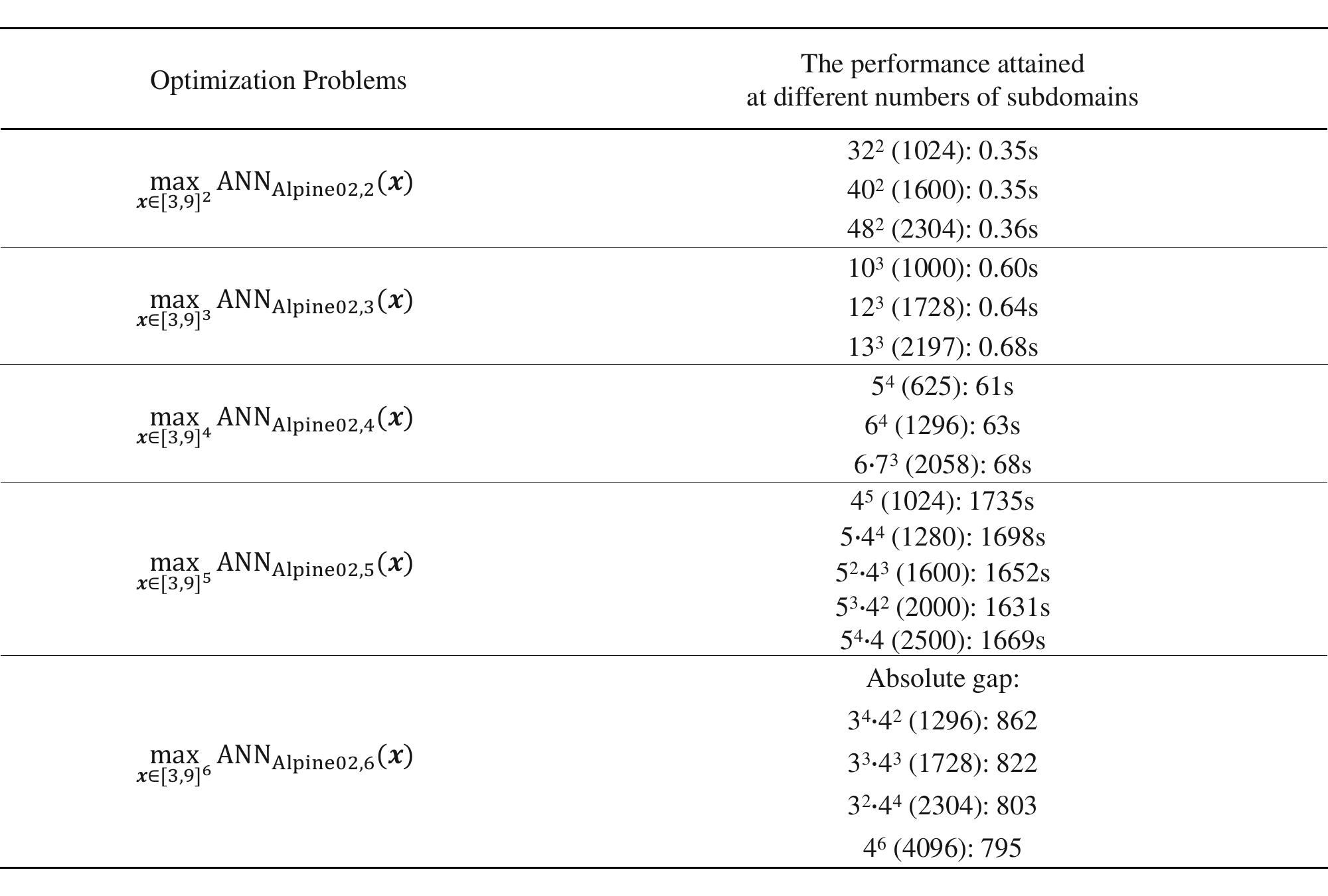}
		\label{tab:Alpine02ANN_Scaling}
	\end{table}
	
	\begin{table}[htbp]
		\captionsetup{position=above}
		\caption{\footnotesize The optimization problems that are used to compare the proposed heterogeneous spatial B\&B framework against MAiNGO's default solver that leverages McCormick relaxations for lower bounding (cf. Section~\ref{para:bound_tightening_test}). For details on the Kinetic ODE parameter estimation problem, refer to \cite{mitsos2009mccormick}. The last column shows the performances attained at different numbers of subdomains (regarding wall-clock time). Though in Table~\ref{tab:Alpine02ANN_Scaling} using 4096 subdomains results in the best absolute gap in the six-dimensional case, for simpler problems shown here, the gain on bound tightness originating from more subdomains cannot compensate for the overhead resulting from more serialization.}
		\centering
		\includegraphics[width=\linewidth, trim=0cm 1.5cm 0cm 1.5cm, clip]{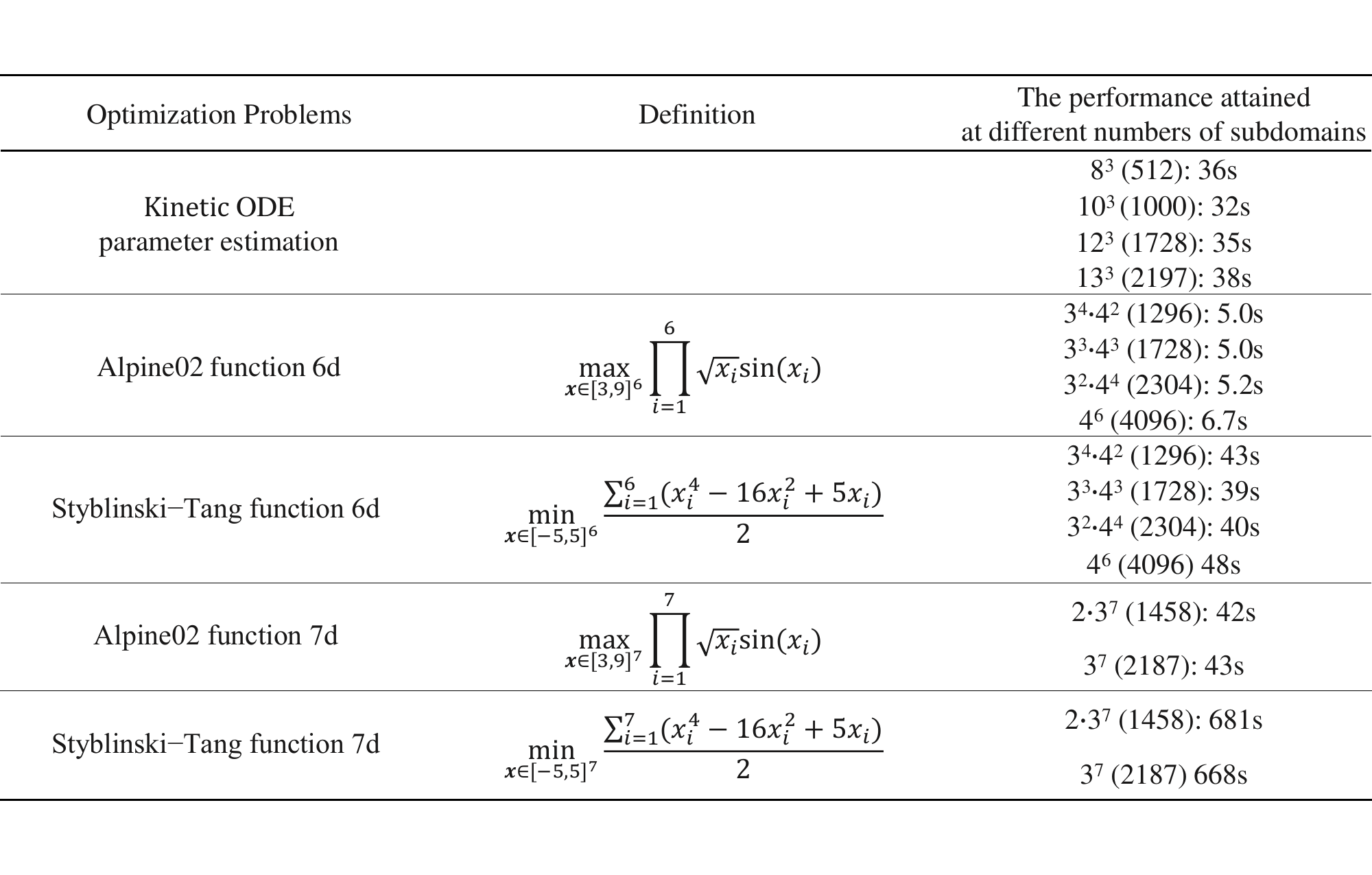}
		\label{tab:ComparisonProblems_details}
	\end{table}
	
	\newpage
	\section*{Statements and Declarations}
	
	\textbf{Author Contributions}
	All authors contributed to the conception and design of this study. Hongzhen Zhang, Tim Kerkenhoff, Manuel Dahmen and Dominik Bongartz formulated the methodology of Subdomain Lower Bounding. Tests, data collection and analysis were performed by Hongzhen Zhang and Tim Kerkenhoff. The first idea of the manuscript was completed by Tim Kerkenhoff and the first draft was written by Hongzhen Zhang. Neil Kichler contributed to the explanation on the interval library and the AD library in the manuscript. Dominik Bongartz offered consistent supports on the manuscript structuring. All authors contributed to the revision of previous versions of the manuscript. The main software of Subdomain Lower Bounding was implemented by Hongzhen Zhang and Tim Kerkenhoff. Supporting libraries used in the CUDA Graph implementation, including the refined interval library, the AD library, and the converter from DAG to CUDA Graph, are implemented by Neil Kichler. All authors read and approved the final manuscript.
	\\\\
	\noindent\textbf{Code Availability}
	The source code of this study has been made publicly available at \url{https://git.rwth-aachen.de/avt-svt/public/maingo}.
	\\\\
	\noindent\textbf{Data Availability}
	The datasets generated during and/or analysed during the current study are available from the corresponding author on reasonable request.
	\\\\
	\noindent\textbf{Competing Interests}
	The authors have no relevant financial or non-financial interests to disclose.

\end{document}